\documentclass[11pt]{amsart}
\usepackage{amsmath,amssymb,amsthm,mathrsfs}
\usepackage{geometry,booktabs}
\usepackage[colorlinks=true,linkcolor=blue,citecolor=blue,urlcolor=blue]{hyperref}
\numberwithin{equation}{section}
\newtheorem{theorem}{Theorem}[section]
\newtheorem{proposition}[theorem]{Proposition}
\newtheorem{lemma}[theorem]{Lemma}
\newtheorem{corollary}[theorem]{Corollary}

\newcommand{\kk}{\Bbbk}
\newcommand{\C}{\mathbb{C}}
\newcommand{\Z}{\mathbb{Z}}
\newcommand{\Q}{\mathbb{Q}}
\newcommand{\N}{\mathbb{N}}
\newcommand{\id}{\operatorname{id}}
\newcommand{\B}{\mathfrak{B}}
\newcommand{\Pp}{\mathcal{P}}
\newcommand{\YD}{\mathcal{YD}}
\newcommand{\HH}{\operatorname{HH}}
\newcommand{\Tor}{\operatorname{Tor}}
\newcommand{\Supp}{\operatorname{Supp}}
\newcommand{\gr}{\operatorname{gr}}
\newcommand{\ord}{\operatorname{ord}}
\newcommand{\Hilb}{\operatorname{Hilb}}
\newcommand{\Rel}{\operatorname{Rel}}
\newcommand{\supp}{\operatorname{supp}}
\newcommand{\End}{\operatorname{End}}

\title[Pre-Nichols algebras and pointed Hopf algebras]{Pre-Nichols algebras and generation of pointed Hopf algebras in positive characteristic}
\author{Rongchuan Xiong}
\address{Department of Mathematics, Changzhou University, Changzhou 213164, China}
\email{rcxiong@foxmail.com}

\makeatletter
\@namedef{subjclassname@2020}{\textup{2020} Mathematics Subject Classification}
\makeatother

\begin{document}
	
	\begin{abstract}
		We study when a finite-dimensional pre-Nichols algebra $R$ of $V$ over an algebraically closed field of characteristic $p>0$, with $p\nmid\dim R$, is necessarily isomorphic to the Nichols algebra $\B(V)$. A trace identity gives a rigidity criterion involving the twist on minimal defining relations. Using integral forms and Hochschild cohomology, we establish uniform rigidity for Yetter--Drinfeld modules of bounded dimension over fixed finite abelian groups, and obtain related results for certain modules over nonabelian groups.
		
		For every finite group $G$ of odd order, outside finitely many characteristics depending only on $G$, every finite-dimensional pointed Hopf algebra $H$ with $G(H)\simeq G$ and dimension prime to the characteristic is generated by group-like and skew-primitive elements. If $G$ is abelian, then, outside a possibly larger finite set of characteristics, every such $H$ is a Hopf $2$-cocycle deformation of its associated graded Hopf algebra.
	\end{abstract}

	\maketitle

	\section{Introduction}

	A basic problem in the theory of finite-dimensional pointed Hopf
	algebras is whether they are generated by their group-like and
	skew-primitive elements. Andruskiewitsch and Schneider formulated
	this as a conjecture in characteristic zero
	\cite[Conjecture~5.7]{AS}. Angiono proved the conjecture when the
	group of group-like elements is abelian
	\cite{Angiono}. More recently, Andruskiewitsch,
	Heckenberger, and Vendramin proved it for pointed Hopf algebras
	of odd dimension in characteristic zero. They also showed that
	these algebras are cocycle deformations of bosonizations of
	Nichols algebras of diagonal type \cite{AHV}.
	
	We consider the generation problem in positive characteristic.
	For a finite-dimensional pointed Hopf algebra $H$ over a field
	$\kk$ of characteristic $p>0$, write $G=G(H)$.
	Its coradical filtration gives
	\[
	\gr H\simeq Q\#\kk G,
	\]
	where $Q$ is the diagram of $H$, a connected graded braided Hopf
	algebra in the Yetter--Drinfeld category over $\kk G$.
	It satisfies $\Pp(Q)=Q(1)$. Set $V:=Q(1)$.
	The braided Hopf subalgebra of $Q$ generated by $V$ is
	the Nichols algebra $\B(V)$ of $V$. Thus $H$ is generated by group-like and
	skew-primitive elements if and only if $Q=\B(V)$.
	
	The graded left dual $Q^*$ is generated by $V^*$
	and is therefore a pre-Nichols algebra of $V^*$.
	This reduces the generation problem to the study of finite-dimensional
	pre-Nichols algebras. Such an algebra has the form
	\[
	R=T(V)/I\twoheadrightarrow
	\B(V)=T(V)/J(V),\qquad I\subseteq J(V),
	\]
	where $I$ is a homogeneous braided Hopf ideal with
	$I\cap V=0$. The question is whether $R$ can be
	finite-dimensional when $I\ne J(V)$.
	
	In positive characteristic, proper finite-dimensional
	pre-Nichols algebras occur even in rank one.
	Let $\kk x$ have self-braiding parameter $q$ of order $N>1$.
	Then, for every $a\ge1$,
	\[
	R=\kk[x]/(x^{Np^a})
	\twoheadrightarrow
	\kk[x]/(x^N)=\B(\kk x)
	\]
	is a proper finite-dimensional pre-Nichols algebra.
	Here $x^N$ is primitive and $\dim R=Np^a$.
	This example leads us to ask whether a finite-dimensional
	pre-Nichols algebra must equal its Nichols algebra when
	its dimension is not divisible by the characteristic.

	Our first result gives a uniform rigidity theorem for
	Yetter--Drinfeld modules over finite abelian groups.
	
	\begin{theorem}\label{thm:uniform}
		Let $G$ be a finite abelian group and $d\geq 1$.
		There is a finite set $S(G,d)$ of rational primes,
		containing the prime divisors of $|G|$, with the
		following property.
		
		Suppose $\kk$ is algebraically closed of characteristic
		$p\notin S(G,d)$ and
		$V\in{}_{\kk G}^{\kk G}\YD$ satisfies $\dim V\leq d$.
		Let $R$ be a finite-dimensional pre-Nichols algebra
		of $V$ with $p\nmid\dim R$. Then $R\simeq\B(V)$.
	\end{theorem}

	The exceptional set $S(G,d)$ can be controlled using integral forms and torsion in third Hochschild cohomology. The universal coefficient sequence relates this torsion to certain components
	of the minimal relation spaces in the residue fibers.
	Since only finitely many braidings arise for fixed
	$G$ and $d$, this gives a finite set of exceptional
	characteristics.

	We also obtain rigidity results for nonabelian
	groups. One such result,
	Corollary~\ref{cor:centerless-rigidity},
	applies to simple Yetter--Drinfeld modules.
	For $X=\bigoplus_{g\in G}X_g$, write
	\[
	\supp X=\{g\in G:X_g\ne0\}.
	\]
	Suppose $G$ is a finite group with $Z(G)=1$,
	$\kk$ is algebraically closed with $|G|1_\kk\ne0$,
	and $V\in{}_{\kk G}^{\kk G}\YD$ is simple with
	\[
	\langle\supp V\rangle=G,
	\qquad \dim V>4.
	\]
	Then every finite-dimensional pre-Nichols
	algebra $R$ of $V$ satisfies
	\[
	(\dim R)1_\kk\ne0
	\quad\Longrightarrow\quad
	R\simeq\B(V).
	\]
	Thus, in characteristic $p>0$, the conclusion
	holds whenever $p\nmid |G|$ and $p\nmid\dim R$.
	In characteristic zero, it holds for every
	finite-dimensional $R$.
	
	For example, consider the transposition and
	$4$-cycle modules over $S_4$ with constant
	cocycle $-1$. Both are simple of dimension
	six, with one-dimensional homogeneous
	components and support generating $S_4$.
	Their Nichols algebras have dimension $576$.
	They therefore provide finite-dimensional
	examples in every characteristic $p>3$.
	
	Theorem~\ref{thm:structural-uniform}, proved
	in Section~\ref{sec:structural}, extends both
	the abelian theorem and the preceding
	nonabelian result. It applies to direct sums
	of one-dimensional modules supported in central
	degrees and modules of dimension greater than four
	that are absolutely simple over the subgroups
	generated by their supports. Writing
	$H=\langle\supp V\rangle$ and
	$H_i=\langle\supp M_i\rangle$ for the larger summands $M_i$,
	the additional condition is that the image of
	$Z(H)\cap H_i$ in $H/[H,H]$ has odd order. The proof uses
	Theorem~\ref{thm:uniform} for the subalgebra
	generated by the one-dimensional summands
	and rank-two restrictions for the remaining
	relations.
	
	For groups of odd order, these structural
	conditions also yield a uniform result without
	a prescribed bound on $\dim V$.
	Indeed, suppose $G$ is a finite group of odd order
	and $\kk$ is algebraically closed of characteristic
	$p>0$ with $p\nmid |G|$.
	If $\B(V)$ is finite-dimensional and
	$p\nmid\dim\B(V)$, then $V$ has the
	decomposition appearing in
	Theorem~\ref{thm:structural-uniform}, with
	\[
	\dim V\leq 2|G|^2.
	\]
	Consequently, there is a finite exceptional
	set $S_{\mathrm{odd}}(G)$ depending only
	on $G$, outside which the same rigidity
	conclusion holds; see Theorem~\ref{thm:odd-general}.
	For abelian $G$, this also follows from
	Gra\~na's multiplicity bound and change
	of characteristic, as explained in
	Corollary~\ref{cor:odd-uniform}.
	
	Applied to the graded dual of the diagram, these rigidity results give
	the generation assertion in the following theorem.
	
	\begin{theorem}\label{thm:intro-cocycle}
		Let $G$ be a finite group of odd order.
		There is a finite set $S_{\mathrm{odd}}(G)$
		of rational primes such that every finite-dimensional
		pointed Hopf algebra $H$ with $G(H)\simeq G$
		over an algebraically closed field $\kk$
		of characteristic $p\notin S_{\mathrm{odd}}(G)$,
		and with $p\nmid\dim H$, is generated by its
		group-like and skew-primitive elements.
		
		If $G$ is abelian, there is a finite set
		$S_{\mathrm{def}}(G)\supseteq S_{\mathrm{odd}}(G)$
		such that, under the same hypotheses with
		$p\notin S_{\mathrm{def}}(G)$,
		\[
		H\simeq(\gr H)_\sigma
		\]
		for some Hopf $2$-cocycle $\sigma$.
	\end{theorem}

	The generation assertion also implies that pointed finite tensor
	categories admitting fiber functors are tensor-generated by objects
	of composition length at most two under the corresponding
	characteristic and dimension restrictions; see
	Corollary~\ref{cor:pointed-categories}.
	The characteristic-zero
	deformation result is due to Angiono and
	Garc\'ia Iglesias \cite[Corollary~1.2]{AGI}.

	We briefly explain the main ideas of the proofs.
	
	The starting point is a general observation about
	defining relations of Nichols algebras; see, for
	example, \cite{Angiono,XiongSubquandle}.
	Suppose that $R\twoheadrightarrow\B(V)$ first
	fails to be an isomorphism in degree $n$. Then
	\[
	0\ne E_n=
	\ker\bigl(R(n)\to\B(V)(n)\bigr)
	\subseteq\mathcal P(R).
	\]
	
	We use a commutator identity to exclude
	certain primitive elements. Let $0\ne r\in R(n)_g$ be primitive,
	with $n>0$, and suppose that $g\cdot r=r$.
	An operator $T$, defined
	using the coproduct and a linear functional,
	satisfies
	\[
	[T,L_r]=\id_R,
	\]
	where $L_r$ denotes left multiplication
	by $r$. For finite-dimensional $R$,
	taking traces gives
	\[
	(\dim R)1_\kk=0.
	\]
	Thus $p\mid\dim R$ in characteristic $p>0$.
	Related $q$-Weyl relations, involving left
	multiplication and opposite braided
	derivatives on Nichols algebras, appear
	in \cite[Proposition~3.2 and
	Remark~3.9]{XiongSquareFree}.
	
	To determine when this identity applies,
	consider the twist on a Yetter--Drinfeld
	module $X=\bigoplus_{g\in G}X_g$:
	\[
	\Theta_X(x_g)=g\cdot x_g,
	\qquad x_g\in X_g.
	\]
	We show that $\Theta_R^{-1}S_R^2$ is an
	algebra endomorphism. Consequently, if
	$R$ is generated by $V=R(1)$ and
	$\Theta_V=q\,\id_V$, then
	\[
	\Theta_R(r)=q^n r,
	\qquad r\in\mathcal P(R)(n).
	\]
	
	Set
	\[
	\Rel(V)=
	\frac{J(V)}
	{T(V)_+J(V)+J(V)T(V)_+}.
	\]
	The first nonzero kernel $E_n$ is a
	quotient of $\Rel(V)(n)$. Hence, if
	$\Theta_{\Rel(V)}=\id$, every element
	of $E_n$ is fixed by the twist.
	The trace identity therefore gives
	\[
	\Theta_{\Rel(V)}=\id,\quad
	\dim R<\infty
	\quad\Longrightarrow\quad
	R=\B(V)\ \text{or}\ p\mid\dim R
	\]
	in characteristic $p>0$.
	In characteristic zero, it gives
	$R=\B(V)$; see
	Theorem~\ref{thm:twist-rigidity}.
	The condition holds, for example, when
	$\Theta_V=-\id_V$ and $\B(V)$ is
	quadratic.
	
	To pass between characteristics, we use
	integral forms of finite-dimensional
	Nichols algebras. Fix a number field $K$
	and a finite-dimensional $B_K=\B(V_K)$.
	After a suitable localization $A$ of
	the ring of integers of $K$, there is a
	finite free graded braided Hopf algebra
	$B_A$ whose residue fibers are again
	Nichols algebras.
	
	For $B=\B(V)$, Hochschild cohomology with trivial
	coefficients describes the minimal relation space:
	\[
	\HH^2(B,\kk)\simeq\Rel(V)^*.
	\]
	In particular, $E_n^*$ embeds into
	$\HH^2_\kk(B,\kk)_n$.
	For the localization $D=A_{\mathfrak p}$ at a
	nonzero prime ideal $\mathfrak p$ and a residue field
	extension $\kk/(A/\mathfrak p)$, the universal
	coefficient sequence gives
	\[
	0\longrightarrow {}
	\HH_D^2(B_D,D)\otimes_D\kk
	\longrightarrow {}
	\HH_\kk^2(B_\kk,\kk)
	\longrightarrow {}
	\Tor_1^D\bigl(
	\HH_D^3(B_D,D),\kk\bigr)
	\longrightarrow 0.
	\]
	Thus new second cohomology classes
	after reduction are governed by torsion
	in third Hochschild cohomology.
	If the minimal relations in the generic
	fiber are fixed by the twist, this gives
	the rigidity criterion of
	Theorem~\ref{thm:general-arithmetic}
	outside finitely many characteristics depending on the
	chosen integral form.
	
	For diagonal braidings, Angiono's
	presentation \cite{Angiono} determines
	the weights of the minimal relations.
	The reflection theory of
	Heckenberger--Schneider \cite{HS}
	allows us to compare finite-dimensional
	Nichols algebras of dimension prime to
	$p$ across characteristics.
	Together with divisibility of $\dim R$ by the dimensions
	of Nichols algebras of primitive Yetter--Drinfeld submodules,
	these results prove
	Theorem~\ref{thm:uniform}.
	
	For nonabelian groups, we use the rank-two restrictions of
	Heckenberger--Vendramin \cite{HVRankTwo}.
	The primitive subspace $V\oplus E_n$
	gives rise to a finite-dimensional
	Nichols algebra, so these restrictions
	apply to its simple constituents.
	
	The proof of
	Theorem~\ref{thm:structural-uniform}
	uses these restrictions to reduce
	$E_n$ to one-dimensional simple
	constituents supported in central
	group degrees. Relations among the
	one-dimensional summands are excluded
	by Theorem~\ref{thm:uniform}.
	For the larger simple summands,
	the mixed braided commutators have noncentral
	group degrees, so only relations within the
	individual summands can contribute.
	Equivariance and the group grading of tensor products
	force the corresponding self-braiding scalars to satisfy
	$q^2=1$. The hypothesis on central
	elements in
	Theorem~\ref{thm:structural-uniform},
	together with the trace identity,
	excludes these remaining constituents.
	
	The cocycle deformation result follows by specialization
	through the finite collection of
	integral forms, as established in
	Appendix~\ref{app:cocycle}.

	The paper is organized as follows. In
	Section~\ref{sec:general}, we study primitive elements
	and the twist. Section~\ref{sec:integral}
	deals with integral forms and Hochschild cohomology.
	Section~\ref{sec:rigidity} proves a general reduction criterion
	and then treats Nichols algebras of diagonal type.
	Section~\ref{sec:structural} treats relations over finite
	groups and proves the structural and odd-order rigidity results.
	Section~\ref{sec:pointed} applies these results to
	finite-dimensional pointed Hopf algebras.

	\section{Primitive elements and defining relations}\label{sec:general}
	
	Let $\kk$ be a field and $G$ a group. The category
	${}_{\kk G}^{\kk G}\mathcal{YD}$ of left-left
	Yetter--Drinfeld modules over $\kk G$
	\cite[Section~1.2]{AS} has as objects the $G$-graded
	$\kk G$-modules $X$ satisfying
	\[
	X=\bigoplus_{g\in G}X_g,\qquad
	a\cdot X_g\subseteq X_{aga^{-1}}.
	\]
	It is a braided monoidal category, with braiding given by
	\[
	c(x\otimes y)=g\cdot y\otimes x,
	\qquad x\in X_g.
	\]
	Algebras and coalgebras in this category are defined by structure maps
	that are morphisms in ${}_{\kk G}^{\kk G}\mathcal{YD}$.
	A braided bialgebra is a bialgebra in this braided monoidal category.
	Its coproduct is an algebra morphism with respect to the braided tensor
	product, whose multiplication is given by
	\[
	(a\otimes b)(u\otimes v)
	=a(g\cdot u)\otimes bv,\qquad b\in R_g.
	\]
	We also use a nonnegative integer grading that is preserved by all
	structure maps. We say that $R$ is connected if $R(0)=\kk1$, and set
	\[
	\Pp(R)=\{r:\Delta(r)=r\otimes1+1\otimes r\}.
	\]
	
	Fix a finite-dimensional object $V$ in
	${}_{\kk G}^{\kk G}\mathcal{YD}$.
	The tensor algebra $T(V)$ is a braided Hopf algebra in which the
	elements of $V$ are primitive. Let $J=J(V)$ be the largest homogeneous
	braided Hopf ideal whose components in degrees zero and one are zero.
	The quotient
	\[
	B=\B(V)=T(V)/J
	\]
	is the Nichols algebra of $V$, characterized by being generated in
	degree one and satisfying $\Pp(B)=V$; see \cite{Nichols} and
	\cite[Section~2]{AS}. A pre-Nichols algebra is a quotient
	\[
	R=T(V)/I,\qquad I\subseteq J,
	\]
	where $I$ is a homogeneous braided Hopf ideal. Let $T_+=\bigoplus_{n\ge1}T^n(V)$
	denote the augmentation ideal of $T(V)$. Then
	\[
	\Rel(V)=J/(T_+J+JT_+)
	\]
	is the space of minimal homogeneous defining relations.

	Suppose that $\pi:R\twoheadrightarrow B$ is first non-injective in
	degree $n$, and set $E_n=\ker\bigl(\pi|_{R(n)}\bigr)$. For $r\in E_n$,
	\[
	\overline\Delta(r)\in
	\bigoplus_{i=1}^{n-1}R(i)\otimes R(n-i),\qquad
	(\pi\otimes\pi)\overline\Delta(r)=0,
	\]
	where $\overline\Delta(r)=\Delta(r)-r\otimes1-1\otimes r$. Since $\pi$
	is injective in every degree $<n$, the map $\pi\otimes\pi$ is
	injective on the displayed direct sum. Thus $\overline\Delta(r)=0$, and
	\begin{equation}\label{eq:first-relation-object}
		E_n\subseteq\Pp(R),\qquad \Rel(V)(n)\twoheadrightarrow E_n.
	\end{equation}
	The second assertion holds because each element of
	$(T_+J+JT_+)(n)$ involves a relation of degree less than $n$ and
	therefore vanishes in $R$.
	
	Consider the quotient
	\[
	A_n=T(V)/\langle J(j):j<n\rangle .
	\]
	The ideal defining $A_n$ is a Hopf ideal: the reduced coproduct of
	each of its generators lies in the sum of tensor products in which
	one factor lies in $J(j)$ for some $j<n$.
	Every element of $J(n)$ has primitive image in $A_n$, since its
	reduced coproduct vanishes modulo the lower-degree relations.
	Conversely, let $x\in A_n$ be a primitive element of degree $n$.
	Its image in $B$ is primitive of degree $n$; since $n\ge2$ and
	$\Pp(B)=V$, this image is zero. Hence
	\[
	\Pp(A_n)(n)=\Rel(V)(n);
	\]
	see also \cite[Lemma~2.1 and Theorem~2.2]{XiongSubquandle}.
	
	If $h$ is the top degree of a finite-dimensional $B$, then
	$J(h+1)=T^{h+1}(V)$. Therefore the ideal generated by $J(j)$ for
	$j\le h+1$ contains every higher tensor degree, and hence equals
	$J$. Consequently
	\begin{equation}\label{eq:relation-degree-bound}
		R\ne B\quad\Longrightarrow\quad 2\le n\le h+1.
	\end{equation}
	
	Assume that $G$ is finite. Then the Nichols--Zoeller theorem
	\cite{NZ}, applied to bosonizations and, for a quotient, to the
	corresponding dual Hopf algebras, gives
	\begin{equation}\label{eq:bosonized-divisibility}
		U\hookrightarrow R\ \text{or}\ R\twoheadrightarrow U
		\quad\Longrightarrow\quad
		\dim U\mid\dim R
	\end{equation}
	for finite-dimensional braided Hopf algebras in
	${}_{\kk G}^{\kk G}\mathcal{YD}$.
	
	\begin{lemma}\label{lem:primitive-factor}
		Assume $G$ is finite, and $R$ is a finite-dimensional connected
		graded braided Hopf algebra in ${}_{\kk G}^{\kk G}\mathcal{YD}$.
		Suppose $W\subseteq\Pp(R)$ is a graded subobject, i.e.\
		$W=\bigoplus_{n\ge1}(W\cap R(n))$, in
		${}_{\kk G}^{\kk G}\mathcal{YD}$. Then
		\[
		\dim\B(W)\mid\dim R.
		\]
	\end{lemma}
	
	\begin{proof}
		Denote by $U\subseteq R$ the braided Hopf subalgebra generated by
		$W$. The filtration
		\[
		F_jU=\sum_{i=0}^j W^i,\qquad W^0=\kk1,
		\]
		is a Hopf filtration in ${}_{\kk G}^{\kk G}\mathcal{YD}$, because $W$
		consists of primitive elements and is stable under the
		Yetter--Drinfeld structure. Since $W$ lies in positive degrees, we have
		$W\cap\kk1=0$, hence $F_1U/F_0U\simeq W$. Therefore $\gr_F U$ is
		a connected graded braided Hopf algebra generated by its degree-one
		part $W$, i.e.\ a pre-Nichols algebra of $W$, and consequently
		there is a surjection $\gr_F U\twoheadrightarrow\B(W)$. By
		\eqref{eq:bosonized-divisibility},
		\[
		\dim\B(W)\mid\dim\gr_F U=\dim U\mid\dim R.
		\]
	\end{proof}
	
	For the first kernel $E_n$ in \eqref{eq:first-relation-object}, the
	subobject $V\oplus E_n\subseteq\Pp(R)$ is a direct sum, since $V$
	sits in degree $1$ while $E_n$ sits in degree $n\ge2$. Because
	$V$ already generates $R$, the filtration by products of
	$V\oplus E_n$ applied in Lemma~\ref{lem:primitive-factor} yields
	\begin{equation}\label{eq:primitive-enlargement}
		\gr_F R\twoheadrightarrow\B(V\oplus E_n),\qquad
		\dim\B(V\oplus E_n)\mid\dim R.
	\end{equation}
	In particular, $E_n$ remains nonzero in filtration degree one, as
	$F_1R/F_0R\simeq V\oplus E_n$. This is the filtration defined in
	\cite[Definition~2.1 and Remark~2.2]{HS}.
	
	For $r\in R$, write $L_r\in\End_\kk(R)$ for left
	multiplication by $r$, so that $L_r(a)=ra$.
	The following operator calculation yields restrictions on homogeneous
	primitive elements in finite-dimensional braided Hopf algebras.
	\begin{lemma}\label{lem:trace}
		Let $G$ be any group and $R$ a finite-dimensional connected graded
		braided Hopf algebra in ${}_{\kk G}^{\kk G}\mathcal{YD}$. Suppose
		\[
		0\ne r\in\Pp(R)\cap R(n)_g,\qquad n>0,\qquad g\cdot r=qr.
		\]
		For $X=L_r$, there is an operator $T$ on $R$ such that
		$TX-q^{-1}XT=\id_R$. If $q=1$, then $(\dim R)1_\kk=0$;
		in particular, in characteristic $p>0$, one has $p\mid\dim R$.
		If $\operatorname{char}\kk=0$, then $q$ has finite order $N>1$,
		$r^N=0$, and $N\mid\dim R$.
	\end{lemma}
	\begin{proof}
		Choose a linear functional $\varphi\in R(n)^*$ with $\varphi(r)=1$
		and extend it by zero in the other degrees. Put
		\[
		U=(\id\otimes\varphi)\Delta,\qquad X=L_r,\qquad G_g(a)=g\cdot a.
		\]
		The braided multiplicativity of $\Delta$ gives
		\[
		\Delta(ra)=
		\sum ra_{(1)}\otimes a_{(2)}
		+\sum(g\cdot a_{(1)})\otimes ra_{(2)}.
		\]
		In the second sum, $\varphi(ra_{(2)})$ vanishes unless $a_{(2)}$
		has degree zero. The corresponding coproduct component
		is $a\otimes1$. Thus $UX=XU+G_g$. Since $G_gX=qXG_g$,
		the operator $T=G_g^{-1}U$ satisfies
		\[
		TX-q^{-1}XT=\id_R.
		\]
		When $q=1$, taking ordinary traces gives
		$0=\operatorname{tr}(\id_R)=(\dim R)1_\kk$.
		
		Now suppose $\operatorname{char}\kk=0$.
		Since $X$ raises degree by $n>0$, there is
		$M\ge1$ such that $X^M=0$ and $X^{M-1}\ne0$.
		Iteration of the operator identity gives
		\[
		TX^M-q^{-M}X^MT=[M]_{q^{-1}}X^{M-1},
		\qquad [m]_a=\sum_{j=0}^{m-1}a^j.
		\]
		Hence $[M]_{q^{-1}}=0$. The case $q=1$ has been excluded,
		so $q$ has finite order $N>1$, with $N\mid M$.
		The quantum binomial formula gives
		$\Delta(r^N)=r^N\otimes1+1\otimes r^N$.
		Moreover, $r^N$ has group degree $g^N$ and
		$g^N\cdot r^N=q^{N^2}r^N=r^N$.
		If $r^N\ne0$, the trace argument just proved applies to $r^N$,
		a contradiction. Thus $r^N=0$ and $M=N$.
		Since $[j]_{q^{-1}}\ne0$ for $1\le j<N$,
		\cite[Theorem~3.1]{XiongSquareFree}, with parameter $q^{-1}$,
		makes $R$ free over $\kk[t]/(t^N)$ via $t\cdot a=Xa$.
		Consequently $N\mid\dim R$.
	\end{proof}

	The argument also applies if there is a graded algebra
	automorphism $\sigma$ of $R$ satisfying
	\[
	c(r\otimes a)=\sigma(a)\otimes r,\qquad
	\sigma(r)=qr \quad (a\in R).
	\]
	Indeed, one replaces $G_g$ by $\sigma$ and uses $\sigma^N$
	for $r^N$.
	
	In particular, let $V$ be a braided vector space of diagonal
	type with basis $(x_i)_{i\in I}$ and
	\[
	c(x_i\otimes x_j)=q_{ij}x_j\otimes x_i.
	\]
	For a graded braided Hopf quotient $R$ of $T(V)$,
	suppose $r\in\Pp(R)$ is nonzero and represented by tensors
	of multidegree $\alpha=(\alpha_i)_{i\in I}$. Then
	\[
	\sigma(x_i)=\left(\prod_{j\in I}q_{ji}^{\alpha_j}\right)x_i,
	\qquad
	\sigma(r)=q_{\alpha\alpha}r,
	\quad
	q_{\alpha\alpha}
	=\prod_{i,j\in I}q_{ij}^{\alpha_i\alpha_j}.
	\]
	The braiding and its inverse show that $\sigma$ preserves
	the defining ideal, so it induces an automorphism of $R$
	satisfying the preceding identities.

	In ${}_{G}^{G}\mathcal{YD}$, we use the twist
	\begin{equation}\label{eq:canonical-twist}
		\Theta_X(x_g)=g\cdot x_g,\qquad
		\Theta_{X\otimes Y}
		=c_{Y,X}c_{X,Y}(\Theta_X\otimes\Theta_Y).
	\end{equation}
	The first formula defines natural automorphisms, and the second is
	the balancing identity; see \cite[Definition~8.10.1]{EGNO}.
	The tensor identity follows on
	$x_g\otimes y_h$: both sides are
	$(gh)\cdot x_g\otimes(gh)\cdot y_h$.
	
	\begin{lemma}\label{lem:antipode-twist}
		For a graded braided Hopf algebra $R$ in ${}_{\kk G}^{\kk G}\mathcal{YD}$,
		$\Theta_R^{-1}S_R^2$ is an algebra endomorphism. If $R$ is generated
		by $V=R(1)$ and $\Theta_V=q\,\id_V$, then
		\[
		\Theta_R|_{R(n)}=q^nS_R^2|_{R(n)},\qquad
		\Theta_R(r)=q^n r\quad(r\in\Pp(R)(n)).
		\]
	\end{lemma}
	\begin{proof}
		Write $m$ for multiplication and $c=c_{R,R}$. The braided
		antipode identity is $S_Rm=mc(S_R\otimes S_R)$.
		Applying it twice and using naturality of $c$ gives the first
		identity below; \eqref{eq:canonical-twist} gives the second:
		\[
		S_R^2m=mc^2(S_R^2\otimes S_R^2),\qquad
		\Theta_Rm=mc^2(\Theta_R\otimes\Theta_R).
		\]
		Thus, with $F=\Theta_R^{-1}S_R^2$,
		\[
		Fm=\Theta_R^{-1}mc^2(S_R^2\otimes S_R^2)
		=m(\Theta_R^{-1}\otimes\Theta_R^{-1})
		(S_R^2\otimes S_R^2)
		=m(F\otimes F),\qquad F(1)=1.
		\]
		On $V$, the identity $S_R(v)=-v$ gives $F=q^{-1}\id$.
		Generation in degree one gives $F|_{R(n)}=q^{-n}\id$, proving the
		first formula. Every primitive satisfies $S_R^2(r)=r$, which gives
		the second.
	\end{proof}
	
	\begin{theorem}\label{thm:twist-rigidity}
		Let $G$ be any group, and suppose that
		$\Theta_{\Rel(V)}=\id$.
		Then every finite-dimensional pre-Nichols algebra $R$ of $V$
		with $(\dim R)1_\kk\ne0$ equals $\B(V)$.
		In particular, this holds in characteristic zero,
		and in characteristic $p>0$ when $p\nmid\dim R$.
	\end{theorem}
	\begin{proof}
		Suppose $R\ne\B(V)$, and choose the least degree $n$ at which the
		kernel is nonzero. Then
		\eqref{eq:first-relation-object} and naturality give
		\[
		0\ne E_n\subseteq\Pp(R),\qquad
		\Theta_{E_n}=\id.
		\]
		Choose $0\ne r\in(E_n)_g$. Then
		$g\cdot r=\Theta_R(r)=r$, so
		Lemma~\ref{lem:trace} gives
		$(\dim R)1_\kk=0$, a contradiction.
	\end{proof}
	
	\begin{corollary}\label{cor:degree-rigidity}
		Let $\Theta_V=q\,\id_V$.
		Suppose $J(V)$ is generated by homogeneous relations in degrees $n$
		satisfying $q^n=1$.
		Then every finite-dimensional pre-Nichols algebra $R$
		of $V$ with $(\dim R)1_\kk\ne0$ equals $\B(V)$.
		In particular, this holds if $\Theta_V=-\id_V$
		and $\B(V)$ is quadratic.
	\end{corollary}
	\begin{proof}
		The classes of the chosen homogeneous generators of $J(V)$ span
		$\Rel(V)$, because multiplication by $T_+$ vanishes in the quotient.
		Their degrees satisfy $q^n=1$, so $\Rel(V)(n)=0$ when $q^n\ne1$.
		In the quotient $A_n$ above, $\Rel(V)(n)$ consists of degree-$n$
		primitive elements. Lemma~\ref{lem:antipode-twist} therefore gives
		\[
		\Theta_{\Rel(V)(n)}=q^n\id.
		\]
		Thus $\Theta_{\Rel(V)}=\id$, and
		Theorem~\ref{thm:twist-rigidity} applies.
	\end{proof}
	
	For finite groups, we also use the following consequence of the
	rank-two classification. If $M,N$ are absolutely simple over a
	nonabelian group generated by their combined supports, then
	\begin{equation}\label{eq:rank-two-bound}
		\begin{gathered}
			\dim\B(M\oplus N)<\infty,\quad c_{N,M}c_{M,N}\ne\id\\
			\Longrightarrow\quad
			\min\{\dim M,\dim N\}\le2,\qquad
			\max\{\dim M,\dim N\}\le4.
		\end{gathered}
	\end{equation}
	This holds in every characteristic by
	\cite[Theorem~2.1 and Table~1]{HVRankTwo} and also \cite[Proposition~3.4]{AHV}.
	
	\begin{proposition}\label{prop:relative-rigidity}
		Suppose $G$ is finite, $\kk$ is algebraically closed with
		$|G|1_\kk\ne0$, and $R$ is a finite-dimensional pre-Nichols
		algebra of $V$. Choose submodules $D_j\subseteq V$ such that the
		subalgebra generated by each $D_j$ in $R$ is $\B(D_j)$. Put
		\[
		Q=\Rel(V)\Big/\sum_j\operatorname{im}
		\bigl(\Rel(D_j)\longrightarrow\Rel(V)\bigr).
		\]
		Suppose each simple constituent $E$ of $Q$ satisfies one of the
		following conditions:
		\begin{enumerate}
			\item $\Theta_E=\id$;
			\item for some $L\le G$ there are absolutely simple submodules
			$M\subseteq\bigoplus_{g\in L}V_g$ and
			$N\subseteq\bigoplus_{g\in L}E_g$ over $L$, whose combined
			supports generate the nonabelian group $L$, with
			\[
			c_{N,M}c_{M,N}\ne\id,\qquad
			\min\{\dim M,\dim N\}>2\ \text{or}\
			\max\{\dim M,\dim N\}>4.
			\]
		\end{enumerate}
		If $(\dim R)1_\kk\ne0$, then $R=\B(V)$.
	\end{proposition}
	\begin{proof}
		The inclusion of tensor algebras induces the map on minimal
		relations, since $J(D_j)=J(V)\cap T(D_j)$ by quantum symmetrizers.
		Suppose $R\ne\B(V)$, and let $E_n$ be the first nonzero homogeneous
		kernel. The quotient map from $\Rel(V)(n)$ vanishes on the images of
		$\Rel(D_j)(n)$ and factors through
		$Q$. Semisimplicity lets us choose a simple constituent $E$ both
		as a submodule of $E_n$ and as a quotient of $Q$.
		In case (1), a nonzero group-homogeneous primitive in $E$
		contradicts Lemma~\ref{lem:trace}. In case (2), the corresponding
		copy of $M\oplus N$ in $V\oplus E_n$ has a finite-dimensional Nichols algebra
		by \eqref{eq:primitive-enlargement}, contradicting
		\eqref{eq:rank-two-bound} over $L$.
	\end{proof}

	\section{Integral forms and Hochschild cohomology}\label{sec:integral}
	
	Fix a finite group $G$, a number field $K$, and
	$V_K\in{}_{G}^{G}\YD$ with finite-dimensional Nichols algebra
	$B_K$. Write $h$ for its top degree and $\mathcal O_K$ for the ring of
	integers of $K$. An integral form permits comparison with residue
	characteristics using a single finite algebra.
	
	\begin{proposition}\label{prop:general-model}
		There is $0\ne a\in\mathcal O_K$ such that, over
		$A=\mathcal O_K[1/(a|G|)]$, the algebra $B_K$ has a finite free
		graded braided Hopf algebra model $B_A$, generated by $V_A=B_A(1)$, with
		\[
		B_A\otimes_AK\simeq B_K,\qquad
		B_A\otimes_A\kappa\simeq\B(V_A\otimes_A\kappa)
		\]
		for every residue field $\kappa$ of $A$.
	\end{proposition}
	\begin{proof}
		Choose a $K$-basis of $B_K$ homogeneous for the integer and group
		gradings, including $1$ and a basis of $V_K$. The matrices of
		multiplication, coproduct, antipode, counit and the finitely many
		group actions have finitely many coefficients. Clearing their
		denominators gives a free lattice $B_A\subset B_K$ preserved by
		all structure maps. Their identities hold over $A$ because the
		lattice and its tensor powers embed in the corresponding spaces
		over $K$.
		
		For $1\leq n\leq h$, generation gives a generically surjective map
		\[
		\mu_n:V_A^{\otimes n}\longrightarrow B_A(n).
		\]
		For $2\leq n\leq h$, the condition $\Pp(B_K)=B_K(1)$ gives a
		generically injective map
		\[
		\overline\Delta_n:B_A(n)\longrightarrow
		\bigoplus_{i=1}^{n-1}B_A(i)\otimes_A B_A(n-i).
		\]
		For each map of positive rank, invert one nonzero maximal minor.
		There are only finitely many such maps. After this localization,
		every residue map $\mu_n$ is surjective and every
		$\overline\Delta_n$ is injective. Since $B_A(n)=0$ for $n>h$,
		each residue fiber satisfies
		\[
		B_\kappa=\langle B_\kappa(1)\rangle,\qquad
		\Pp(B_\kappa)=B_\kappa(1).
		\]
		The characterization of Nichols algebras proves the assertion.
		The product of the elements inverted can be incorporated into $a$.
	\end{proof}
	
	Write $B_L=B_A\otimes_A L$ and $V_L=V_A\otimes_A L$.
	The cohomological arguments below apply as well to a finite
	projective model generated in degree one whose generic and residue
	fibers are Nichols algebras.
	
	For an augmented algebra $S$ over a commutative ring $L$, we use
	Hochschild cohomology with coefficients in $L$, regarded as a
	trivial $S$-bimodule via the augmentation. With
	$S_+=\ker\epsilon$, the normalized Hochschild cochain complex
	\cite[Section~9.1 and Exercise~8.6.4]{Weibel} is
	\begin{equation}\label{eq:bar}
		\begin{split}
			\mathcal C_L^j(S)&=\operatorname{Hom}_L(S_+^{\otimes_Lj},L),\\
			(df)(a_1,\ldots,a_{j+1})
			&=\sum_{i=1}^j(-1)^i
			f(a_1,\ldots,a_i a_{i+1},\ldots,a_{j+1}).
		\end{split}
	\end{equation}
	The internal degree is the total input degree. When $S$ is a
	finite projective Yetter--Drinfeld algebra, the cochains have the
	left dual structure in the same category \cite[Section~1.2]{AS}:
	\begin{equation}\label{eq:dual-convention}
		\deg_G(f|_{X_g})=g^{-1},\qquad
		(a\cdot f)(x)=f(a^{-1}\cdot x),\qquad
		(X\otimes Y)^*\simeq Y^*\otimes X^*.
	\end{equation}
	Each summand of the bar differential is the dual of a
	Yetter--Drinfeld multiplication map. It therefore preserves the
	action and group grading. Direct evaluation gives
	\[
	(\Theta_{X^*}f)(x_g)=f(g\cdot x_g),\qquad
	\Theta_{X^*}=\Theta_X^{\mathsf t},
	\]
	where the superscript denotes the linear transpose.
	
	\begin{proposition}\label{prop:relation-module}\label{prop:first-class}
		Let $B=\B(V)$ over a field $\kk$. In each internal degree
		$n\geq2$ there is a natural isomorphism of Yetter--Drinfeld modules
		\[
		\HH_\kk^2(B,\kk)_n\simeq\Rel(V)(n)^*.
		\]
		If $B$ is finite-dimensional, these assemble into a graded isomorphism
		\begin{equation}\label{eq:minimal-relations}
			\HH_\kk^2(B,\kk)\simeq\Rel(V)^*.
		\end{equation}
		For a proper pre-Nichols quotient $R\twoheadrightarrow B$,
		let $n$ and $E_n$ be as in \eqref{eq:first-relation-object}. Then
		\begin{equation}\label{eq:first-class}
			E_n^*\lhook\joinrel\longrightarrow\HH_\kk^2(B,\kk)_n.
		\end{equation}
		If $B$ has top degree $h$, then $2\leq n\leq h+1$.
	\end{proposition}
	\begin{proof}
		Decomposing tensor words at their last letter gives
		\[
		Z=\ker(B\otimes V\to B_+)\simeq J/JV,\qquad JV=JT_+.
		\]
		The exact sequences
		\[
		0\to B_+\to B\to\kk\to0,\qquad
		0\to Z\to B\otimes V\to B_+\to0
		\]
		give, by dimension shifting,
		\[
		\operatorname{Tor}_2^B(\kk,\kk)
		\simeq\kk\otimes_B Z
		\simeq J/(T_+J+JT_+).
		\]
		Indeed, the map $\kk\otimes_BZ\to\kk\otimes_B(B\otimes V)=V$
		is zero, because $J$ has degree at least two. The normalized bar
		resolution identifies each internal component of $\HH^2$ with
		the dual of the corresponding component of this
		$\operatorname{Tor}_2$. If $B$ is finite-dimensional, this gives
		\eqref{eq:minimal-relations};
		compare \cite[Section~4.1]{MW}. The exact sequences and the bar
		resolution are natural for the action and grading, giving the
		asserted Yetter--Drinfeld compatibility.
		
		Dualizing \eqref{eq:first-relation-object} gives
		\eqref{eq:first-class}; \eqref{eq:relation-degree-bound} gives the
		degree bound. For injections defined by reduced coproducts in coalgebra
		Hochschild cohomology, see e.g. \cite{StefanVO,WZZ,XiongP2q}.
	\end{proof}
	
	Fix a nonzero prime ideal $\mathfrak p\subset A$, and put
	\[
	D=A_{\mathfrak p},\qquad \kappa=A/\mathfrak p,\qquad
	p=\operatorname{char}\kappa.
	\]
	The ring $D$ is a discrete valuation ring. The terms of
	$\mathcal C_D^\bullet(B_D)$ are finite free, and their
	base changes are the normalized complexes over $K$ and over every
	field extension $\kk/\kappa$.
	
	\begin{proposition}\label{prop:base-change}
		The module $\HH_D^2(B_D,D)$ is torsion-free. For every field
		extension $\kk/\kappa$, there is a natural exact sequence
		\begin{equation}\label{eq:UCT}
			\begin{split}
				0\longrightarrow \HH_D^2(B_D,D)\otimes_D\kk
				&\longrightarrow \HH_\kk^2(B_\kk,\kk)\\
				&\longrightarrow
				\Tor_1^D\bigl(\HH_D^3(B_D,D),\kk\bigr)
				\longrightarrow0.
			\end{split}
		\end{equation}
		It preserves internal degrees and Yetter--Drinfeld structures.
	\end{proposition}
	\begin{proof}
		For $n\geq2$, multiplication is a surjection
		\[
		\bigoplus_{i=1}^{n-1}B_D(i)\otimes_D B_D(n-i)
		\twoheadrightarrow B_D(n).
		\]
		It splits over $D$ because its target is free. Its dual is, up to
		sign, $d^1$ in internal degree $n$. Thus $\operatorname{im}d^1$
		is a direct summand of $\mathcal C_D^2$, and
		\[
		\HH_D^2(B_D,D)=\ker d^2/\operatorname{im}d^1
		\subseteq\mathcal C_D^2/\operatorname{im}d^1
		\]
		is torsion-free. Internal degrees below two have no degree-two
		cochains.
		
		For a uniformizer $\varpi$, multiplication by $\varpi$ gives an
		exact sequence of complexes
		\[
		0\to\mathcal C_D^\bullet
		\xrightarrow{\ \varpi\ }\mathcal C_D^\bullet
		\to\mathcal C_\kappa^\bullet\to0.
		\]
		The long exact sequence gives
		\[
		0\to\HH_D^2/\varpi\HH_D^2
		\to\HH_\kappa^2\to\HH_D^3[\varpi]\to0.
		\]
		The identification
		\[
		\HH_D^3[\varpi]\simeq\Tor_1^D(\HH_D^3,\kappa)
		\]
		and extension of scalars to $\kk$ give \eqref{eq:UCT}.
		This universal coefficient sequence for tensor product is also
		obtained by reindexing the K\"unneth formula in
		\cite[Theorem~3.6.1]{Weibel}. All chain maps preserve the action,
		group grading and internal degree; so do the connecting map and
		the scalar extension. This proves the stated compatibilities.
	\end{proof}
	
	Reduction can decrease the ranks of the bar differentials and
	produce additional classes in second cohomology. The quotient by
	reductions of integral classes is measured by the $\Tor_1$ term
	in \eqref{eq:UCT}. On a direct summand where the generic $\HH^2$
	vanishes, torsion-freeness makes the first term zero, so the
	entire residue component is $\Tor_1$ of integral $\HH^3$.
	Such additional classes do not by themselves imply the existence
	of a proper finite-dimensional pre-Nichols algebra.
	
	For finite $G$, the twist defines natural cohomological summands.
	Write $m=\exp G$. Since $m$ is
	invertible in $A$ and $\Theta_X^m=\id$ on every Yetter--Drinfeld
	module over $AG$,
	\begin{equation}\label{eq:twist-projector}
		e=\frac1m\sum_{j=0}^{m-1}\Theta^j
	\end{equation}
	is the projector onto $\ker(\Theta-\id)$. Naturality makes it a
	cochain projector. Its defining formula commutes with localization
	and scalar extension, so it acts compatibly on
	\eqref{eq:UCT}.
	
	\section{Reduction modulo prime ideals}\label{sec:rigidity}
	
	\subsection{Relations fixed by the twist}
	
	Retain the integral form $B_A$ from Section~\ref{sec:integral}, with
	generic fiber $B_K$ of top degree $h$. The condition on the twist and
	base change in Hochschild cohomology give the following result.
	
	\begin{theorem}\label{thm:general-arithmetic}
		Suppose $\Theta_{\Rel(V_K)}=\id$. For every nonzero prime
		ideal $\mathfrak p\subset A$, every field extension $\kk/\kappa$, and
		every $n\geq2$,
		\begin{equation}\label{eq:general-torsion}
			(1-e)\Rel(V_\kk)(n)^*
			\simeq
			\Tor_1^D\bigl((1-e)\HH_D^3(B_D,D)_n,\kk\bigr).
		\end{equation}
		If a proper finite-dimensional pre-Nichols algebra $R$ of $V_\kk$
		has $p\nmid\dim R$, its first nonzero homogeneous kernel satisfies
		\[
		0\ne E_n^*\hookrightarrow
		\Tor_1^D\bigl((1-e)\HH_D^3(B_D,D)_n,\kk\bigr)
		\quad\text{for some }2\leq n\leq h+1.
		\]
		There is consequently a finite set $S(V_A)$ of rational primes
		such that every residue extension of characteristic
		$p\notin S(V_A)$ satisfies
		\[
		\dim R<\infty,\quad p\nmid\dim R
		\quad\Longrightarrow\quad R=\B(V_\kk).
		\]
	\end{theorem}
	\begin{proof}
		By Proposition~\ref{prop:relation-module}, the twist acts trivially
		on $\HH_K^2(B_K,K)$. Proposition~\ref{prop:base-change} shows that
		$\HH_D^2(B_D,D)$ is torsion-free. Its direct summand on which $e$
		acts by zero has zero generic fiber and therefore vanishes:
		\[
		(1-e)\HH_D^2(B_D,D)=0.
		\]
		Apply $1-e$ to \eqref{eq:UCT} and use
		\eqref{eq:minimal-relations} in the residue field to obtain
		\eqref{eq:general-torsion}.
		
		For $E_n=\ker(R(n)\to B_\kk(n))$ in the least degree $n$ where it is nonzero,
		\eqref{eq:first-relation-object} gives $E_n\subseteq\Pp(R)$.
		If $eE_n\ne0$, choose a nonzero group-homogeneous element
		$r\in(eE_n)_g$. Then $g\cdot r=r$, contradicting
		Lemma~\ref{lem:trace}. Thus $eE_n=0$. The equivariance of
		\eqref{eq:first-class}, together with
		$\Theta_{E_n^*}=\Theta_{E_n}^{\mathsf t}$, places its image in
		$(1-e)\HH_\kk^2$. This proves the injection.
		The bound on $n$ is \eqref{eq:relation-degree-bound}.
		
		The module
		\[
		M=(1-e)\bigoplus_{n=2}^{h+1}\HH_A^3(B_A,A)_n
		\]
		is finitely generated over the Dedekind ring $A$. A nonzero ideal
		annihilates its torsion, so only the finitely many prime ideals dividing
		that ideal can occur in the support of the torsion submodule. For this fixed model, let
		$S(V_A)$ consist of the residue characteristics of prime ideals
		of $\mathcal O_K$ removed by localization, together with those below
		$\Supp_A\operatorname{tors}_A M$. Outside this set every
		$\Tor_1$ group in the injection vanishes, so a proper $R$ is
		impossible.
	\end{proof}
	
	More precisely, for $2\leq n\leq h+1$, \eqref{eq:general-torsion}
	gives
	\begin{equation}\label{eq:twist-support}
		\Supp_A\operatorname{tors}_A
		\bigl((1-e)\HH_A^3(B_A,A)_n\bigr)
		=
		\{\mathfrak p\in\operatorname{Spec}A\setminus\{0\}:
		(1-e)\Rel(V_{A/\mathfrak p})(n)\ne0\}.
	\end{equation}
	Thus the torsion support is exactly the locus where the minimal
	relations have a nonzero component in the complement of the fixed
	subspace of the twist.
	
	\begin{corollary}\label{cor:integral-presentation}
		Let $J_A=\ker(T_A(V_A)\to B_A)$. If $J_A$ is generated as a
		two-sided ideal by a homogeneous Yetter--Drinfeld submodule
		$M\subseteq J_A$ with $\Theta_M=\id$, then
		\[
		\operatorname{tors}_A
		\left((1-e)\bigoplus_{n=2}^{h+1}\HH_A^3(B_A,A)_n\right)=0.
		\]
		At every residue extension, a finite-dimensional pre-Nichols
		algebra of dimension prime to the characteristic is Nichols.
	\end{corollary}
	\begin{proof}
		Since $B_A$ is projective, the presentation remains exact after
		reduction. The image of $M\otimes_A\kk$ generates $J(V_\kk)$,
		so it surjects onto $\Rel(V_\kk)$ and makes its twist the identity.
		Theorem~\ref{thm:twist-rigidity} gives the assertion. Equation
		\eqref{eq:twist-support}, applied at every nonzero prime ideal, gives
		the torsion vanishing.
	\end{proof}
	
	For braidings of diagonal type, we restrict the weights of
	primitive elements representing defining relations in
	finite-dimensional pre-Nichols algebras. Dimension divisibility
	and reflection theory restrict the weights of the first nonzero kernel; the twist
	need not act as the identity on the entire space of minimal
	defining relations.
	
	\subsection{Diagonal braidings and change of characteristic}
	\label{sec:diagonal}\label{sec:diagonal-model}
	
	Suppose $G$ is finite abelian of exponent $m$, and $F$ is a field
	containing all $m$th roots of unity, with $|G|1_F\ne0$. Write
	$\widehat G=\operatorname{Hom}(G,F^\times)$. Every finite-dimensional
	Yetter--Drinfeld module is a direct sum of one-dimensional Yetter--Drinfeld modules
	$F_\gamma$, where $\gamma=(g,\chi)\in G\times\widehat G$:
	\[
	\delta(x)=g\otimes x,\qquad a\cdot x=\chi(a)x.
	\]
	For a weight basis $x_1,\ldots,x_\theta$,
	\begin{equation}\label{eq:braiding}
		c(x_i\otimes x_j)=q_{ij}x_j\otimes x_i,\qquad
		q_{ij}=\chi_j(g_i)\in\mu_m(F).
	\end{equation}
	This is a braiding of diagonal type. For a weight $\gamma=(g,\chi)$,
	put $q_{\gamma\gamma}=\chi(g)$. A one-dimensional braided vector space with parameter $q$
	has Nichols algebra
	\begin{equation}\label{eq:rank-one}
		\B(Fx)=
		\begin{cases}
			F[x]/(x^N),&q\ne1,\quad N=\ord(q),\\
			F[x]/(x^p),&q=1,\quad\operatorname{char}F=p>0,\\
			F[x],&q=1,\quad\operatorname{char}F=0.
		\end{cases}
	\end{equation}
	
	Let $\kk$ be algebraically closed of characteristic $p>0$, with
	$p\nmid|G|$, and fix an isomorphism
	$\mu_m(\kk)\simeq\mu_m(\C)$. The primes dividing $m$ are exactly
	those dividing $|G|$, so all these roots retain their orders.
	For $q=(q_{ij})$ over $\kk$, denote its corresponding complex
	matrix by $\widehat q$.
	
	\begin{proposition}\label{prop:transfer}
		Under these assumptions,
		\[
		\dim\B_\kk(q)<\infty,\quad p\nmid\dim\B_\kk(q)
		\quad\Longleftrightarrow\quad
		\dim\B_\C(\widehat q)<\infty.
		\]
		In this case the multigraded Hilbert series agree.
	\end{proposition}
	\begin{proof}
		We use \cite[Theorems~6.7, 7.2 and 7.3]{HS}. Their root
		factorization is available because every tensor power of $V$ is
		semisimple in the present Yetter--Drinfeld category.
		
		For a diagonal braiding, the $i$th reflection is determined by
		$a_{ii}=2$ and
		\begin{equation}\label{eq:cartan}
			-a_{ij}=
			\min\{a\geq0:(a+1)_{q_{ii}}
			(1-q_{ii}^a q_{ij}q_{ji})=0\},\quad i\ne j,
		\end{equation}
		where $(n)_q=\sum_{r=0}^{n-1}q^r$. It replaces the simple degrees by
		\[
		s_i(\alpha_j)=\alpha_j-a_{ij}\alpha_i.
		\]
		The reflected matrix is obtained by evaluating the bicharacter
		$q_{\alpha,\beta}=\prod_{j,l}q_{jl}^{\alpha_j\beta_l}$ on these
		new simple degrees. If $q_{ii}\ne1$, then
		\[
		(a+1)_{q_{ii}}=0\quad\Longleftrightarrow\quad q_{ii}^{a+1}=1.
		\]
		Hence both \eqref{eq:cartan} and the reflected matrix depend only
		on multiplicative relations among roots of unity.
		
		Suppose first that $\B_\kk(q)$ is finite-dimensional of dimension
		prime to $p$. Its iterated braided adjoint actions vanish in sufficiently high degree, so
		all reflections are defined. They preserve dimension by
		\cite[Theorem~6.7(3)]{HS}. In every reflected matrix $q'$, a
		one-dimensional submodule in degree one consists of primitive elements.
		Lemma~\ref{lem:primitive-factor} gives
		\[
		\dim\B_\kk(q'_{ii})\mid\dim\B_\kk(q')=\dim\B_\kk(q).
		\]
		By \eqref{eq:rank-one}, $q'_{ii}\ne1$. Induction on a reflection
		sequence now transfers all Cartan matrices and reflections to
		$\C$. Their integer reflection maps agree, so the Weyl groupoids
		agree and are finite. Every reflected
		one-dimensional complex module has a finite-dimensional Nichols algebra. By
		\cite[Theorem~7.3]{HS}, $\dim\B_\C(\widehat q)<\infty$.
		
		Conversely, finite-dimensionality over $\C$ gives all reflections
		and a finite Weyl groupoid. Every reflected diagonal coefficient
		differs from $1$, by the characteristic-zero case of
		\eqref{eq:rank-one}. Thus \eqref{eq:cartan} transfers all the
		reflection data to $\kk$, where the reflected one-dimensional modules again have
		finite-dimensional Nichols algebras. The same criterion gives
		$\dim\B_\kk(q)<\infty$.
		
		By \cite[Theorem~7.2]{HS}, each module in the root factorization
		is isomorphic to one of the reflected one-dimensional modules.
		Write $\Delta_+\subset\N_0^\theta$ for the
		common positive root set, and put
		\[
		q_{\beta\beta}=\prod_{i,j}q_{ij}^{\beta_i\beta_j},
		\qquad N_\beta=\ord(q_{\beta\beta}).
		\]
		All $q_{\beta\beta}$ differ from $1$, and their orders are
		preserved. The root factorization gives
		\begin{equation}\label{eq:root-product}
			\Hilb_{\B}(\boldsymbol t)=
			\prod_{\beta\in\Delta_+}
			\left(\sum_{r=0}^{N_\beta-1}\boldsymbol t^{r\beta}\right).
		\end{equation}
		The Hilbert series agree, and
		$\dim\B_\kk(q)=\prod_{\beta\in\Delta_+}N_\beta$ is prime to $p$
		because each $N_\beta$ divides $m$.
	\end{proof}
	
	For diagonal data the quantum shuffle construction gives an
	integral form at every prime ideal of residue characteristic not dividing the group order.
	Set
	\begin{equation}\label{eq:coefficient-ring}
		\mathcal O=\Z[\zeta_m,1/m],\qquad K=\Q(\zeta_m),
	\end{equation}
	and choose a direct sum $V_{\mathcal O}$ of one-dimensional Yetter--Drinfeld modules over
	$\mathcal O G$ with finite-dimensional $B_K=\B(V_K)$.
	The character idempotents split the weight decomposition over
	$\mathcal O$. On $V_{\mathcal O}^{\otimes n}$ let
	$\mathfrak S_n=\sum_{\sigma\in S_n}c_\sigma$ be the quantum
	symmetrizer, using the positive braid lifts of permutations, and set
	\begin{equation}\label{eq:integral-model}
		J_{\mathcal O}(n)=\ker\mathfrak S_n,\qquad
		B_{\mathcal O}
		=T_{\mathcal O}(V_{\mathcal O})/
		\bigoplus_{n\geq0}J_{\mathcal O}(n).
	\end{equation}
	
	\begin{proposition}\label{prop:model}
		The algebra $B_{\mathcal O}$ is a finite projective connected
		graded braided Hopf algebra, generated by $V_{\mathcal O}$, with
		generic fiber $B_K$. For every nonzero prime ideal
		$\mathfrak p\subset\mathcal O$ and residue field
		$\kappa=\mathcal O/\mathfrak p$,
		\[
		B_{\mathcal O}\otimes_{\mathcal O}\kappa\simeq\B(V_\kappa),
		\qquad
		\Hilb_{\B(V_\kappa)}(\boldsymbol t)
		=\Hilb_{B_K}(\boldsymbol t).
		\]
		Moreover $\operatorname{char}\kappa\nmid\dim\B(V_\kappa)$.
	\end{proposition}
	\begin{proof}
		Fix $\mathfrak p$, put $D=\mathcal O_{\mathfrak p}$, and let
		$p=\operatorname{char}\kappa$. Reduction identifies the $m$th roots
		of unity in $K$ and $\kappa$. Over a field $L$, the degree-$n$
		quantum shuffle map has image $\B(V_L)(n)$ and is
		$\mathfrak S_n\otimes L$ \cite{Rosso}. This description holds in every
		characteristic: the iterated shuffle product is the quantum symmetrizer,
		and the deconcatenation coproduct has no primitive elements in degrees
		greater than one. Hence its subalgebra generated by $V_L$ is $\B(V_L)$.
		Proposition~\ref{prop:transfer} therefore gives, in every multidegree,
		\[
		\operatorname{rank}_K\mathfrak S_n
		=\operatorname{rank}_\kappa(\mathfrak S_n\otimes\kappa).
		\]
		Write $b$ for the common rank on one such component. A $b$-minor
		is a unit in $D$, whereas all larger minors vanish. Row and column
		operations over $D$ give
		\[
		\mathfrak S_n\otimes D\sim
		\begin{pmatrix}I_b&0\\0&0\end{pmatrix}.
		\]
		Thus its kernel, image and cokernel are free, and its kernel
		commutes with reduction. At every nonzero prime ideal this proves
		\[
		B_{\mathcal O}(n)\simeq\operatorname{im}\mathfrak S_n
		\text{ is projective},\qquad
		J_{\mathcal O}(n)\otimes\kappa
		=\ker(\mathfrak S_n\otimes\kappa).
		\]
		
		The quantum shuffle map is a braided Hopf algebra morphism \cite{Rosso}.
		Its image is locally a direct summand in each tensor degree, as
		just proved; its tensor-square inclusion is therefore injective.
		Compatibility with coproducts makes the kernel of the quantum shuffle map a coideal as well
		as an ideal. Hence \eqref{eq:integral-model} has the induced
		braided Hopf structure and the stated generic and residue fibers.
		
		If $h$ is the top degree of $B_K$, then
		$\mathfrak S_n\otimes K=0$ for $n>h$. The tensor lattice is
		torsion-free, so $\mathfrak S_n=0$ in those degrees.
		Thus $B_{\mathcal O}$ is finite projective.
		The Hilbert series and the final assertion follow
		from Proposition~\ref{prop:transfer}.
	\end{proof}
	
	\subsection{Weights of relations and exceptional primes}
	
	Throughout this subsection, retain the cyclotomic integral form above,
	with top degree $h$. We use the left dual convention
	\eqref{eq:dual-convention}. For a weight $\gamma=(g,\chi)$, its dual
	weight is therefore
	\begin{equation}\label{eq:dual-weight}
		\gamma^*=(g^{-1},\chi^{-1})\in G\times\widehat G.
	\end{equation}
	All cohomology weight subscripts below refer to this convention.
	
	\begin{proposition}\label{prop:generic-weights}
		Let $K$ be a characteristic-zero splitting field for $G$, with
		$B=\B(V)$ finite-dimensional. If $\gamma=(g,\chi)$ satisfies
		\[
		\chi(g)\ne1,\qquad\dim\B(V\oplus K_\gamma)<\infty,
		\]
		then
		\[
		\HH_K^2(B,K)_{n,\gamma^*}=0\qquad(n\geq2).
		\]
	\end{proposition}
	\begin{proof}
		Choose homogeneous generators from
		\cite[Theorem~3.1]{Angiono} whose classes form a basis of
		$\Rel(V)$. They generate $J(V)$ by induction on degree.
		For each such generator $\rho$ of weight $\gamma_\rho$,
		Angiono's calculations for primitive relations give
		\begin{equation}\label{eq:relation-weight-exclusion}
			q_{\rho\rho}\ne1
			\quad\Longrightarrow\quad
			\dim\B(V\oplus K_{\gamma_\rho})=\infty.
		\end{equation}
		Appendix~\ref{app:weights} supplies the weight formulas and the
		precise correspondence with \cite[Section~4]{Angiono}.
		It follows that no minimal generator has weight $\gamma$.
		The assertion is now \eqref{eq:minimal-relations}.
		
		For an arbitrary splitting field $K$ of characteristic zero,
		the weight data are defined over a cyclotomic number field
		$K_0\subseteq K$. The homogeneous components of $\B(V)$ and the minimal
		relation space $\Rel(V)$ commute with field extension. Applying the preceding
		calculation over $\C$ and descending to $K_0$ proves the result
		over $K$.
	\end{proof}
	
	Let
	\begin{equation}\label{eq:E}
		\mathcal E(V_K)=
		\{\gamma=(g,\chi)\in G\times\widehat G:
		\chi(g)\ne1,\ \dim\B(V_K\oplus K_\gamma)<\infty\}.
	\end{equation}
	Every weight occurring in $E_n$ for a finite-dimensional pre-Nichols
	algebra of dimension prime to $p$ belongs to this set, after lifting
	the roots of unity to $K$.
	
	\begin{theorem}\label{thm:first-torsion}
		Fix a nonzero prime ideal $\mathfrak p\subset\mathcal O$, and set
		$D=\mathcal O_{\mathfrak p}$, $\kappa=\mathcal O/\mathfrak p$,
		$p=\operatorname{char}\kappa$. For an algebraically
		closed extension $\kk$ of $\kappa$, if $R\twoheadrightarrow B_\kk$ is a proper finite-dimensional
		pre-Nichols algebra with $p\nmid\dim R$, then
		$2\leq n\leq h+1$ for the least degree with nonzero kernel,
		and the dual of that kernel admits a weight-preserving injection
		\begin{equation}\label{eq:torsion-injection}
			E_n^*\lhook\joinrel\longrightarrow
			\bigoplus_{\gamma\in\mathcal E(V_K)}
			\Tor_1^D\bigl(\HH_D^3(B_D,D)_{n,\gamma^*},\kk\bigr).
		\end{equation}
	\end{theorem}
	\begin{proof}
		By \eqref{eq:first-relation-object},
		$E_n\subseteq\Pp(R)$. For a nonzero
		$r\in(E_n)_\gamma$, with $\gamma=(g,\chi)$,
		Lemma~\ref{lem:primitive-factor} gives
		\[
		\dim\B(\kk r)\mid\dim R.
		\]
		Consequently $\chi(g)\ne1$, by \eqref{eq:rank-one}.
		Since $n\geq2$, the submodule $V_\kk\oplus\kk r$ is a direct sum
		of primitives. Applying the same lemma gives
		\[
		\dim\B(V_\kk\oplus\kk r)\mid\dim R.
		\]
		Proposition~\ref{prop:transfer} makes
		$\B(V_K\oplus K_\gamma)$ finite-dimensional, so
		$\gamma\in\mathcal E(V_K)$.
		
		The generic $\gamma^*$-component of $\HH^2$ vanishes by
		Proposition~\ref{prop:generic-weights}. Its integral counterpart
		is torsion-free by Proposition~\ref{prop:base-change}, hence
		\[
		\HH_D^2(B_D,D)_{n,\gamma^*}=0.
		\]
		Apply the $\gamma^*$ weight projector to \eqref{eq:UCT}. This yields
		\[
		\HH_\kk^2(B_\kk,\kk)_{n,\gamma^*}
		\simeq\Tor_1^D\bigl(\HH_D^3(B_D,D)_{n,\gamma^*},\kk\bigr).
		\]
		Taking the direct sum over the weights of $E_n$ and using
		\eqref{eq:first-class} proves \eqref{eq:torsion-injection}.
	\end{proof}
	
	Define
	\begin{equation}\label{eq:torsion-module}
		\mathcal T(V_{\mathcal O})=
		\operatorname{tors}_{\mathcal O}
		\left(\bigoplus_{n=2}^{h+1}
		\bigoplus_{\gamma\in\mathcal E(V_K)}
		\HH_{\mathcal O}^3(B_{\mathcal O},\mathcal O)_{n,\gamma^*}
		\right).
	\end{equation}
	The finite exceptional set for these fixed weight data is
	\begin{equation}\label{eq:fixed-primes}
		S(V_{\mathcal O})={}\{p:p\mid|G|\}
		{}\cup
		\{\operatorname{char}(\mathcal O/\mathfrak p):
		\mathfrak p\in\Supp_{\mathcal O}
		\mathcal T(V_{\mathcal O})\}.
	\end{equation}
	
	\begin{theorem}\label{thm:arithmetic-main}
		The set $S(V_{\mathcal O})$ is finite. For every reduction
		$V_\kk$ over an algebraically closed field of characteristic
		$p\notin S(V_{\mathcal O})$ and every finite-dimensional
		pre-Nichols algebra $R$ of $V_\kk$,
		\[
		p\nmid\dim R\quad\Longrightarrow\quad R=\B(V_\kk).
		\]
	\end{theorem}
	\begin{proof}
		The normalized cochain terms for $B_{\mathcal O}$ are finite
		projective. Their cohomology, and each weight component of it,
		are finitely generated over the Dedekind domain $\mathcal O$.
		The sum in \eqref{eq:torsion-module} is finite, so its torsion
		has finite support. This proves finiteness of $S(V_{\mathcal O})$.
		
		At a prime ideal above $p\notin S(V_{\mathcal O})$, every local
		cohomology component in \eqref{eq:torsion-injection} is
		torsion-free and its $\Tor_1$ vanishes. A proper $R$ would
		therefore have $E_n^*=0$, contrary to the choice of $n$.
	\end{proof}
	
	\begin{corollary}\label{cor:local-test}
		At a nonzero prime ideal $\mathfrak p\subset\mathcal O$, suppose
		$\HH_D^3(B_D,D)_{n,\gamma^*}$ is torsion-free for every
		$2\leq n\leq h+1$ and $\gamma\in\mathcal E(V_K)$.
		Then every finite-dimensional pre-Nichols algebra of $V_\kk$
		with $p\nmid\dim R$ is Nichols.
	\end{corollary}
	\begin{proof}
		Every $\Tor_1$ term in \eqref{eq:torsion-injection} is zero.
	\end{proof}
	
	\begin{proof}[Proof of Theorem~\ref{thm:uniform}]
		Write $\widehat G=\operatorname{Hom}(G,K^\times)$. For each
		multiset of at most $d$ weights from $G\times\widehat G$, take
		the direct sum of the corresponding one-dimensional
		Yetter--Drinfeld modules over $\mathcal O G$.
		Let $\mathcal V(G,d)$ be the finite collection of these
		$V_{\mathcal O}$ whose generic Nichols algebras are
		finite-dimensional, and set
		\[
		S(G,d)=\{p:p\mid|G|\}\cup
		\bigcup_{V_{\mathcal O}\in\mathcal V(G,d)}S(V_{\mathcal O}).
		\]
		Each set in this union is finite by
		Theorem~\ref{thm:arithmetic-main}, so $S(G,d)$ is finite and
		depends only on $G$ and $d$.
		
		Suppose $V$ and $R$ satisfy the hypotheses of the theorem.
		Choosing a primitive $m$th root of unity in $\kk$ defines a
		homomorphism $\mathcal O\to\kk$, whose kernel is a nonzero
		prime ideal $\mathfrak p$. Reduction identifies the $m$th roots
		of unity in $K$ and $\kk$, so the weight decomposition of $V$
		lifts to a direct sum $V_{\mathcal O}$ of the form chosen above,
		with $V_{\mathcal O}\otimes_{\mathcal O}\kk\simeq V$.
		
		The quotient $R\twoheadrightarrow\B(V)$ and
		\eqref{eq:bosonized-divisibility} give
		$\dim\B(V)\mid\dim R$. Thus $\B(V)$ is finite-dimensional of
		dimension prime to $p$, and Proposition~\ref{prop:transfer}
		makes $\B(V_K)$ finite-dimensional. Hence
		$V_{\mathcal O}\in\mathcal V(G,d)$. Since
		$p\notin S(V_{\mathcal O})$ and $\kk$ is an extension of
		$\mathcal O/\mathfrak p$, Theorem~\ref{thm:arithmetic-main}
		gives $R=\B(V)$.
	\end{proof}
	
	\begin{corollary}\label{cor:odd-uniform}
		For every finite abelian group $G$ of odd order, there is a finite
		set $S(G)$ of rational primes such that, over every algebraically closed field of
		characteristic $p\notin S(G)$, every finite-dimensional
		pre-Nichols algebra over $\kk G$ satisfies
		\[
		p\nmid\dim R\quad\Longrightarrow\quad R=\B(R(1)).
		\]
	\end{corollary}
	\begin{proof}
		By \cite[Proposition~3.2]{GranaLow}, each simple summand of a characteristic-zero object
		with a finite-dimensional Nichols algebra over an odd-order group occurs at
		most twice. In its proof,
		repeated copies give a constant diagonal matrix $q_{ij}=q$.
		Finite-dimensionality excludes $q=1$, and odd group order
		excludes $q=-1$; the bound follows from the resulting Cartan
		matrix. If $G$ is abelian, there are $|G|^2$ one-dimensional Yetter--Drinfeld modules, so
		\[
		\dim V_K\leq2|G|^2.
		\]
		For the $R$ in the statement,
		$\dim\B(R(1))\mid\dim R$. Proposition~\ref{prop:transfer}
		therefore gives the same dimension bound after lifting to
		characteristic zero. Take $S(G)=S(G,2|G|^2)$ and apply
		Theorem~\ref{thm:uniform}.
	\end{proof}
	
	\section{Relations over finite groups}\label{sec:structural}
	
	Throughout this section $\kk$ is algebraically closed, $G$ is
	finite, and $|G|1_\kk\ne0$. Set $H=\langle\supp V\rangle$.
	Every group degree of $R(n)$ is a product of $n$ elements of
	$\supp V$, since $T^n(V)\twoheadrightarrow R(n)$ preserves group
	degrees. We may therefore restrict the action and coaction to $H$;
	the braiding and $J(V)$ are unchanged. This restriction
	does not in general preserve simplicity. Simplicity over a support
	subgroup will be specified whenever it is used.
	
	\subsection{Mixed monodromy and central degrees}
	
	\begin{lemma}\label{lem:local-monodromy}
		Let $M\in{}_{\kk H}^{\kk H}\YD$ be absolutely simple over
		$H_M=\langle\supp M\rangle$, with $\dim M>4$, and suppose $E$ is an
		object over $H$. If $c_{E,M}c_{M,E}\ne\id$, then there are
		$g\in\supp E$, $L=\langle H_M,g\rangle$, and an absolutely simple
		submodule $N\subseteq\bigoplus_{x\in g^L}E_x$ over $L$ such that
		\[
		\langle\supp M,\supp N\rangle=L,\qquad
		c_{N,M}c_{M,N}\ne\id.
		\]
		In particular, if $\B(M\oplus E)$ is finite-dimensional, then its
		mixed monodromy is the identity.
	\end{lemma}
	\begin{proof}
		Choose $g\in\supp E$ such that the mixed monodromy is not the
		identity on $M\otimes E_g$.
		The subspace $\bigoplus_{x\in g^L}E_x$ is a Yetter--Drinfeld module over $L$.
		It is semisimple, and each of its simple summands has that full
		conjugacy orbit as support. The restriction of $M$ to $L$ is
		absolutely simple because it is already so over $H_M\le L$.
		For some summand $N$, the mixed monodromy with $M$ is not the
		identity. Its support contains $g$, so the two supports generate
		$L$. This group is nonabelian: otherwise its simple object $M$
		would be one-dimensional. The finite-dimensional assertion follows
		from \eqref{eq:rank-two-bound}. Here $M\oplus N$ is a braided
		subspace, so its Nichols algebra embeds in $\B(M\oplus E)$.
	\end{proof}
	
	\begin{lemma}\label{lem:central-character}
		Let $M_1,\ldots,M_s$ be submodules of $V$, with $s\ge1$, each
		absolutely simple over $H_i=\langle\supp M_i\rangle$.
		Suppose $H$ is generated by the $H_i$ together with a subset of
		$Z(H)$. If a simple object $E$ over $H$ satisfies
		$c_{E,M_i}c_{M_i,E}=\id$ for every $i$, then
		$E=\kk_g^\chi$ for some $g\in Z(H)$ and character
		$\chi:H\to\kk^\times$.
	\end{lemma}
	\begin{proof}
		For $v_h\in(M_i)_h$ and $e_g\in E_g$, the mixed monodromy is
		\[
		c_{E,M_i}c_{M_i,E}(v_h\otimes e_g)
		=(hgh^{-1})\cdot v_h\otimes h\cdot e_g.
		\]
		Its second group degree is $hgh^{-1}$. The identity therefore
		forces $hg=gh$. Thus $\supp E$ centralizes every $H_i$, and hence
		all of $H$. Simplicity puts $E$ in one central degree $g$ and
		makes its ordinary $H$-representation irreducible.
		
		The action of $g$ on $M_i$ preserves its group grading and commutes
		with the action of $H_i$. Schur's lemma in the Yetter--Drinfeld
		category over $H_i$ makes it a scalar $a_i$. The displayed identity
		now gives
		\[
		h\cdot e=a_i^{-1}e,\qquad h\in\supp M_i,\ e\in E.
		\]
		Every remaining central generator acts scalarly on the irreducible
		ordinary representation $E$. All generators of $H$ consequently
		act by scalars. Irreducibility gives $\dim E=1$, and the existing
		group action defines the character $\chi$.
	\end{proof}
	
	\begin{corollary}\label{cor:centerless-rigidity}
		Suppose $H=\langle\supp V\rangle$ is centerless and $V$ contains
		an absolutely simple submodule $M$ over $H$ with
		$\langle\supp M\rangle=H$ and $\dim M>4$.
		Then every finite-dimensional pre-Nichols algebra $R$ of $V$ with
		$(\dim R)1_\kk\ne0$ is $\B(V)$.
	\end{corollary}
	\begin{proof}
		For a first nonzero kernel $E_n$, \eqref{eq:primitive-enlargement}
		makes $\B(M\oplus E_n)$ finite-dimensional.
		Lemma~\ref{lem:local-monodromy} forces the mixed monodromy to be
		the identity. The group-degree calculation in
		Lemma~\ref{lem:central-character} gives
		\[
		\supp E_n\subseteq H\cap C_H(\langle\supp M\rangle)=Z(H)=1.
		\]
		A nonzero primitive in $E_n$ then contradicts Lemma~\ref{lem:trace}.
	\end{proof}
	
	\subsection{Arithmetic rigidity}
	
	\begin{theorem}\label{thm:structural-uniform}
		Let $G$ be a finite group and $r\ge1$. There is a finite set
		$S_{\mathrm{str}}(G,r)$ of rational primes, containing the prime
		divisors of $|G|$, with the following property. Suppose $\kk$ is
		algebraically closed of characteristic $p\notin S_{\mathrm{str}}(G,r)$
		and $V\in{}_{\kk G}^{\kk G}\YD$ satisfies $\dim V\le r$.
		Set $H=\langle\supp V\rangle$. Suppose that, over $H$,
		\[
		V=D\oplus M_1\oplus\cdots\oplus M_s,
		\]
		where $D$ is a sum of one-dimensional modules of central group
		degrees, and each $M_i$ is absolutely simple over
		$H_i=\langle\supp M_i\rangle$, with $\dim M_i>4$. Suppose also that
		\begin{equation}\label{eq:intro-central-image}
			\left|(Z(H)\cap H_i)[H,H]/[H,H]\right|\ \text{is odd}
			\qquad(1\le i\le s).
		\end{equation}
		Every finite-dimensional pre-Nichols algebra $R$ of $V$ satisfies
		\[
		p\nmid\dim R\quad\Longrightarrow\quad R=\B(V).
		\]
		If $D=0$, the same implication holds in every characteristic
		$p\nmid|G|$. In characteristic zero the same rigidity holds
		without an exceptional set.
	\end{theorem}
	
	For the sets of Theorem~\ref{thm:uniform}, put
	\begin{equation}\label{eq:structural-primes}
		S_{\mathrm{str}}(G,r)=\{p:p\mid|G|\}\cup
		\bigcup_{\substack{C\le G\\C\text{ abelian}}}S(C,r).
	\end{equation}
	Only finitely many subgroups occur, so this set is finite and depends
	only on $G$ and $r$.
	
	\begin{proof}[Proof of Theorem~\ref{thm:structural-uniform}]
		Work over $H=\langle\supp V\rangle$. The subalgebra $R_D$
		generated by $D$ is a pre-Nichols algebra of $D$, and
		\eqref{eq:bosonized-divisibility} gives $\dim R_D\mid\dim R$.
		Its group degrees lie in the abelian subgroup
		$C=\langle\supp D\rangle\le Z(H)$.
		If $D=0$, then $R_D=\B(D)=\kk$.
		For $D\ne0$, restriction to $C$ preserves the braiding.
		In characteristic $p\notin S_{\mathrm{str}}(G,r)$,
		Theorem~\ref{thm:uniform} and \eqref{eq:structural-primes}
		give $R_D=\B(D)$; in characteristic zero this follows from
		\cite[Theorem~4.13]{Angiono}.
		
		Since $\B(V)$ is finite-dimensional, \eqref{eq:rank-two-bound}
		forces identity mixed monodromy between distinct $M_i$ and between
		each $M_i$ and every line of $D$. Both modules remain absolutely
		simple over their combined support subgroup. In particular the
		$H_i$ commute pairwise. The Nichols algebra is the braided tensor
		product
		\begin{equation}\label{eq:structural-product}
			\B(V)\simeq\B(D)\,\underline\otimes\,
			\B(M_1)\,\underline\otimes\cdots\underline\otimes\,\B(M_s).
		\end{equation}
		Its defining ideal is generated by $J(D)$ and the $J(M_i)$
		and the mixed quadratic braided commutators. Indeed their quotient is this
		braided tensor product, generated by $V$ with primitive space $V$.
		Consequently, the images of $\Rel(D)$ and the $\Rel(M_i)$,
		together with the classes of these commutators, span $\Rel(V)$.
		
		If $s=0$, then $R=R_D=\B(V)$. Suppose $s>0$, and let $E$ be a
		simple constituent of $Q=\Rel(V)/\operatorname{im}\Rel(D)$
		that does not satisfy condition (2) of
		Proposition~\ref{prop:relative-rigidity}.
		Lemma~\ref{lem:local-monodromy} makes its mixed monodromy with
		every $M_i$ the identity; otherwise it would give a pair satisfying
		that condition. Lemma~\ref{lem:central-character}
		gives $E=\kk_g^\chi$ with $g\in Z(H)$.
		
		A mixed commutator involving $M_i$ has group degree $xy$, with
		$x\in\supp M_i$ and $y$ centralizing $H_i$.
		If $xy$ were central in $H$, then $x\in Z(H_i)$.
		This is impossible: absolute simplicity over $H_i$ makes $\supp M_i$ a
		nonsingleton conjugacy class generating the nonabelian group
		$H_i$. By semisimplicity, there is a nonzero morphism
		$\Rel(V)(n)\to E$ vanishing on $\Rel(D)(n)$ for some $n$.
		The mixed commutator classes map to zero, so this morphism is
		nonzero on the image of $J(M_i)(n)_g$ for some $i$.
		In particular $g\in Z(H)\cap H_i$.
		
		Write $a$ for the scalar action of $g$ on $M_i$, and
		$b=\chi(x)$ for $x\in\supp M_i$. The latter is constant on
		that conjugacy class. Mixed monodromy gives $ab=1$.
		Equivariance of the nonzero relation map and the group degree
		of its representing tensors give, respectively,
		\begin{equation}\label{eq:relation-character-square}
			\chi(g)=a^n=b^n,\qquad \chi(g)^2=(ab)^n=1.
		\end{equation}
		The order of $\chi(g)$ divides the order of the image of $g$
		in $H/[H,H]$, which is odd by
		\eqref{eq:intro-central-image}. Therefore $\chi(g)=1$ and
		$\Theta_E=\id$. Proposition~\ref{prop:relative-rigidity}, with the
		submodule $D$, gives $R=\B(V)$.
	\end{proof}
	
	Only the odd order of the particular value $\chi(g)$ is used in
	the last step. The group condition
	\eqref{eq:intro-central-image} ensures it uniformly for the
	possible group degrees and characters of relations within each $M_i$.
	
	For a centerless group generated by a conjugacy class $X$ with
	$|X|>4$, every Yetter--Drinfeld module supported on $X$ with
	one-dimensional homogeneous components is absolutely simple:
	a nonzero graded subobject contains one line, and conjugation
	makes it contain them all. Corollary~\ref{cor:centerless-rigidity}
	therefore applies. For $S_4$, both transpositions and $4$-cycles
	generate the group. In the latter case their squares generate
	the normal subgroup of double transpositions, and their images
	in the quotient generate $S_3$. The examples with constant cocycle $-1$ cited
	in the introduction have dimension $576$, prime to every $p>3$.
	For transpositions of $S_n$, $n\ge4$, the same rigidity implication
	holds in characteristic $p\nmid n!$ for finite-dimensional
	pre-Nichols algebras of dimension prime to $p$.
	
	\subsection{Groups of odd order}
	
	\begin{lemma}\label{lem:odd-supports}
		Suppose $G$ has odd order and $\B(V)$ is finite-dimensional of
		dimension nonzero in $\kk$. Then
		\[
		\dim V\le2\sum_{g\in G}(\ord(g)-1)\le2|G|^2,
		\]
		and, over $H=\langle\supp V\rangle$, $V$ has the decomposition
		in Theorem~\ref{thm:structural-uniform}.
	\end{lemma}
	\begin{proof}
		On $V_g$ the action of $g$ is diagonalizable. Each nonzero
		$q$-eigenspace $E$ has constant braiding $q$ times the flip,
		with $q$ of odd order. It is a Yetter--Drinfeld module over
		$\langle g\rangle$, and there is an ordinary Hopf inclusion
		\[
		\B(E)\#\kk\langle g\rangle
		\hookrightarrow\B(V)\#\kk G.
		\]
		The Nichols--Zoeller theorem and $|G|1_\kk\ne0$ show that
		$\dim\B(E)$ is nonzero in $\kk$. The parameter $q=1$ is
		impossible by \eqref{eq:rank-one}. For $q\ne1$, the constant
		matrix lifts to characteristic zero by
		Proposition~\ref{prop:transfer}. Gra\~na's multiplicity bound
		\cite[Proposition~3.2]{GranaLow} then gives $\dim E\le2$:
		it is the repeated one-dimensional object over the odd cyclic
		group of order $\ord(q)$. In characteristic zero apply the
		same bound directly. Summing over the at most $\ord(g)-1$
		nontrivial eigenvalues gives the displayed estimate.
		Decompose $V$ into simple Yetter--Drinfeld modules over $H$.
		The homogeneous component of a summand in group degree $g$ is an irreducible
		representation of $C_H(g)$ on which $g$ acts scalarly, so its
		dimension is at most two. Irreducible degrees of an odd-order
		group are odd, also in nonmodular characteristic; hence this
		homogeneous component is a line. In particular, a simple object supported in a central group degree
		is one-dimensional.
		
		We next show that the supports of distinct simple summands commute.
		If $r,s$ in two such supports do not commute, put
		$L=\langle r,s\rangle$. Their two $L$-conjugacy orbits have odd
		size greater than one, and hence at least three. Choose absolutely simple
		Yetter--Drinfeld submodules over $L$, one in each of the two sums of
		homogeneous components on these orbits. Their dimensions are
		at least three, their supports generate $L$, and their mixed
		monodromy is not the identity because its second group degree changes
		from $s$ to $rsr^{-1}$. Their Nichols algebra embeds in $\B(V)$,
		contradicting \eqref{eq:rank-two-bound}.
		
		Since $H$ is generated by the pairwise commuting support subgroups,
		the support $X$ of a summand $M$ is a single conjugacy class already in
		$H_M=\langle X\rangle$. If $H_M$ is abelian, $X$ is a singleton
		central in $H$ and $M$ is a line, as just proved.
		Otherwise the one-dimensional homogeneous components and the transitive
		$H_M$-action on $X$ show that $M$ is absolutely simple over $H_M$:
		a nonzero graded submodule contains one homogeneous line and hence all of them,
		and the same argument holds after any field extension.
		Its class cannot have size three:
		an odd-order permutation group on three points has no nonidentity
		element fixing a point, whereas each generator $x\in X$ fixes
		$x$. The conjugation action would be trivial, contrary to the
		nonsingleton class $X$. Therefore $\dim M\ge|X|\ge5$.
		Collect the central-degree lines into $D$ and the other summands
		into the $M_i$. Condition~\eqref{eq:intro-central-image} is automatic in
		an odd-order group.
	\end{proof}
	
	\begin{theorem}\label{thm:odd-general}
		For every finite group $G$ of odd order there is a finite set
		$S_{\mathrm{odd}}(G)$ of rational primes such that every finite-dimensional
		pre-Nichols algebra $R$ over $\kk G$, with $\kk$ algebraically
		closed of characteristic $p\notin S_{\mathrm{odd}}(G)$, satisfies
		\[
		p\nmid\dim R\quad\Longrightarrow\quad R=\B(R(1)).
		\]
		One may take
		$S_{\mathrm{odd}}(G)=S_{\mathrm{str}}(G,2|G|^2)$.
		In characteristic zero the same rigidity holds.
	\end{theorem}
	\begin{proof}
		The Nichols quotient has dimension dividing $\dim R$.
		Apply Lemma~\ref{lem:odd-supports} and
		Theorem~\ref{thm:structural-uniform}.
	\end{proof}
	
	\section{Generation of pointed Hopf algebras}\label{sec:pointed}
	
	Let $H$ be a finite-dimensional pointed Hopf algebra with
	$G=G(H)$. The coradical filtration of $H$
	is a Hopf filtration. Write $C=\gr H$ and
	$\pi:C\to\kk G$ for the canonical projection, and put
	\[
	Q=C^{\operatorname{co}\pi}
	=\{x\in C:(\id\otimes\pi)\Delta(x)=x\otimes1\}.
	\]
	The Radford--Majid decomposition gives
	\begin{equation}\label{eq:pointed-diagram}
		\gr H\simeq Q\#\kk G,\qquad V_H=Q(1),
	\end{equation}
	where $Q$ is a connected graded braided Hopf algebra over $\kk G$,
	called the diagram, and $V_H$ is the infinitesimal braiding.
	Moreover, $\Pp(Q)=Q(1)$; see \cite{AS,Radford,Montgomery} for more details.

	The following facts on graded duals of braided Hopf
	algebras are standard; see, e.g.,
	\cite[Lemma~2.4]{AS}.
	\begin{lemma}\label{lem:rigid-hopf-dual}\label{sec:duality}
		Let $Q$ be a finite-dimensional connected graded braided Hopf
		algebra in ${}_{\kk G}^{\kk G}\YD$. Its graded dual
		$Q^*=\bigoplus_n Q(n)^*$, formed using the left dual structure in \eqref{eq:dual-convention}, is a
		braided Hopf algebra in ${}_{\kk G}^{\kk G}\YD$.
		Using $(f\otimes h)(x\otimes y)=f(x)h(y)$, its structure maps are
		\begin{equation}\label{eq:rigid-hopf-dual}
			(fh)(x)=(h\otimes f)\Delta_Q(x),\qquad
			\Delta_{Q^*}(f)(x\otimes y)=f(yx).
		\end{equation}
		If $\Pp(Q)=Q(1)$, then $Q^*$ is a pre-Nichols algebra of $Q(1)^*$.
	\end{lemma}
	\begin{proof}
		The left dual reverses tensor factors. Dualizing the Hopf
		identities therefore gives \eqref{eq:rigid-hopf-dual}, and the
		resulting object is a braided Hopf algebra in the same category.
		The reduced coproduct is dual to multiplication on $Q^*_+$, so
		\[
		\bigl(Q^*_+/(Q^*_+)^2\bigr)^*\simeq\Pp(Q).
		\]
		If $\Pp(Q)=Q(1)$, the indecomposables of $Q^*$ are concentrated
		in degree one; induction on degree then shows that $Q^*$ is generated
		in degree one. Conversely, generation of $Q$ in degree one gives
		$\Pp(Q^*)=Q(1)^*$.
	\end{proof}

	Assume now that $\operatorname{char}\kk=p>0$ and
	$p\nmid\dim H$. The Nichols--Zoeller theorem gives
	$|G|\mid\dim H$, and hence
	\[
	p\nmid|G|,\qquad
	\dim Q^*=\dim Q=\frac{\dim H}{|G|}.
	\]
	Thus the graded dual of the diagram is precisely a finite-dimensional
	pre-Nichols algebra of dimension prime to $p$ to which the
	preceding theorems apply.
	
	\begin{corollary}\label{cor:general-AS}
		Let $H$ be finite-dimensional and pointed in characteristic
		$p>0$, with $p\nmid\dim H$. Set $W:=V_H^*$. If
		$\Theta_{\Rel(W)}=\id$, then $H$ is generated by group-like elements and skew-primitive elements.
		It suffices that $\Theta_{V_H}=q\,\id$ and that $J(W)$ be
		generated in degrees $n$ with $q^n=1$. In odd characteristic,
		$\Theta_{V_H}=-\id$ and quadraticity of $\B(W)$ suffice.
	\end{corollary}
	\begin{proof}
		Lemma~\ref{lem:rigid-hopf-dual} makes $Q^*$ a pre-Nichols algebra
		of $W$. Theorem~\ref{thm:twist-rigidity} gives
		$Q^*=\B(W)$, so $\Pp(Q^*)=W$. Applying the lemma to $Q^*$
		shows that $Q$ is generated by $V_H$. Since $\Pp(Q)=V_H$,
		we obtain $Q=\B(V_H)$. In particular $\gr H$ is generated
		in degrees zero and one.
		
		To lift generation, let $H'$ be the subalgebra generated by $H_1$.
		For $x\in H_n$, express its class in $H_n/H_{n-1}$ as a polynomial
		in degree-zero and degree-one generators of $\gr H$.
		Lift that expression to $a\in H'\cap H_n$. Then
		$x-a\in H_{n-1}$, so induction on $n$ gives $x\in H'$.
		The first term $H_1$ of a pointed coalgebra is spanned by
		group-like and skew-primitive elements, see e.g. \cite[Lemma~2.1]{WZZ}.
		
		The sufficient conditions follow from
		$\Theta_{V_H^*}=\Theta_{V_H}^{\mathsf t}$ and
		Corollary~\ref{cor:degree-rigidity}.
	\end{proof}
	
	If $V_H\simeq V_H^*$, the quadratic condition in this corollary
	can equivalently be checked on $\B(V_H)$.
	
	\begin{corollary}\label{cor:arithmetic-AS}
		Fix a finite group $G$ and $W_K\in{}_{KG}^{KG}\YD$, where $K$ is a number field.
		Suppose $\B(W_K)$ is finite-dimensional and $\Theta_{\Rel(W_K)}=\id$.
		For a fixed integral form $W_A$, outside its finite exceptional
		set, every finite-dimensional pointed Hopf algebra $H$ over a
		residue field extension with
		\[
		G(H)\simeq G,\qquad V_H^*\simeq W_\kk,\qquad p\nmid\dim H
		\]
		is generated by group-like elements and skew-primitive elements.
	\end{corollary}
	\begin{proof}
		Apply Theorem~\ref{thm:general-arithmetic} to the graded dual of the diagram.
		The dualization and filtration argument in
		Corollary~\ref{cor:general-AS} then apply.
	\end{proof}
	
	The diagonal arithmetic theorem gives a uniform generation
	statement.
	\begin{corollary}\label{cor:pointed}
		Fix a finite abelian group $G$ and $d\geq1$. Let $\kk$ be
		algebraically closed of characteristic $p\notin S(G,d)$.
		Every finite-dimensional pointed Hopf algebra $H$ with
		\[
		G(H)\simeq G,\qquad\dim V_H\leq d,\qquad p\nmid\dim H
		\]
		has $\gr H\simeq\B(V_H)\#\kk G$ and is generated by
		group-like and skew-primitive elements. For odd-order $G$, the same
		conclusion holds for $p\notin S(G)$ and arbitrary $\dim V_H$.
	\end{corollary}
	\begin{proof}
		Apply Theorem~\ref{thm:uniform} to $Q^*$, or
		Corollary~\ref{cor:odd-uniform} for odd-order $G$.
		Its degree-one dimension is $\dim V_H$.
		The proof of Corollary~\ref{cor:general-AS} identifies the diagram
		and lifts its generators.
	\end{proof}
	
	\begin{corollary}\label{cor:structural-pointed}
		Fix a finite group $G$ and $r\ge1$. Let $\kk$ be algebraically
		closed of characteristic $p\notin S_{\mathrm{str}}(G,r)$, and
		$H$ a finite-dimensional pointed Hopf algebra with
		$G(H)\simeq G$ and $p\nmid\dim H$. If $\dim V_H\le r$ and
		$V_H$ has the decomposition and satisfies the condition on central elements in
		Theorem~\ref{thm:structural-uniform}, then
		\[
		\gr H\simeq\B(V_H)\#\kk G,\qquad
		H=\langle G(H),\text{skew-primitive elements}\rangle.
		\]
		For every fixed odd-order $G$, the same conclusion holds outside
		$S_{\mathrm{odd}}(G)$ without a bound on $\dim V_H$.
	\end{corollary}
	\begin{proof}
		The structural hypotheses are preserved by the left dual:
		group supports are inverted, their generated subgroups do not
		change, absolute simplicity over those subgroups is preserved, and
		one-dimensional summands in central degrees remain of the same type. Thus
		Theorem~\ref{thm:structural-uniform} applies to the pre-Nichols
		algebra $Q^*$; for an odd-order group, Theorem~\ref{thm:odd-general}
		applies instead. Lemma~\ref{lem:rigid-hopf-dual} and the filtration
		argument in Corollary~\ref{cor:general-AS} give the result.
	\end{proof}
	
	A finite tensor category $\mathcal C$ is pointed if every simple
	object is invertible. Write $G(\mathcal C)$ for the group of
	isomorphism classes of invertible objects and
	$\operatorname{FPdim}(\mathcal C)$ for its Frobenius--Perron
	dimension \cite[Definition~6.1.7]{EGNO}. A fiber functor is an exact
	faithful $\kk$-linear tensor functor to finite-dimensional
	$\kk$-vector spaces.
	
	For a pointed Hopf algebra $H$, the category
	$\operatorname{Comod}_{\mathrm{fd}}(H)$ is pointed, and generation
	of $H$ by group-like and skew-primitive elements is equivalent to
	tensor generation of this category by objects of composition length
	at most two \cite[Proposition~5.11.9]{EGNO}. Here tensor generation
	means that every object is a subquotient of a finite direct sum of
	tensor products of the specified objects
	\cite[Definition~5.11.8]{EGNO}.
	
	\begin{corollary}\label{cor:pointed-categories}
		Let $G$ be a finite group of odd order, and   $\kk$  
		algebraically closed of characteristic
		$p>0$ with $p\notin S_{\mathrm{odd}}(G)$.
		Suppose $\mathcal C$ is a pointed finite tensor category over
		$\kk$ admitting a fiber functor, with
		\[
		G(\mathcal C)\simeq G,\qquad
		p\nmid\operatorname{FPdim}(\mathcal C).
		\]
		Then $\mathcal C$ is tensor-generated by objects of composition
		length at most two.
	\end{corollary}
	\begin{proof}
		Let $F$ be a fiber functor on $\mathcal C$. By
		\cite[Theorem~5.3.12]{EGNO}, $H=\End(F)^*$ is a
		finite-dimensional Hopf algebra and
		$\mathcal C\simeq_\otimes\operatorname{Comod}_{\mathrm{fd}}(H)$.
		Since the simple objects are invertible, their images under $F$
		are one-dimensional. Thus $H$ is pointed and
		$G(H)\simeq G(\mathcal C)\simeq G$. Moreover,
		\[
		\dim H=\operatorname{FPdim}(\mathcal C)
		\]
		by \cite[Proposition~4.5.7 and Example~6.1.9]{EGNO}.
		By Corollary~\ref{cor:structural-pointed}, $H$ is generated by
		group-like and skew-primitive elements, and
		\cite[Proposition~5.11.9]{EGNO} gives the assertion.
	\end{proof}
	
	\begin{corollary}\label{cor:AA-dimension}
		Fix a finite abelian group $G$, $d\geq1$, and a rational prime
		$p\notin S(G,d)$. Suppose $H$ is finite-dimensional and pointed over
		an algebraically closed field of characteristic $p$, with
		$G(H)=G$ and $p\nmid\dim H$. For a Yetter--Drinfeld submodule
		$V'\subseteq V_H$ with $\dim V'\leq d$, choose equivariant
		skew-primitive lifts and let $H'$ be generated by them and $G$.
		Then
		\[
		\dim H'=|G|\dim\B(V').
		\]
		For odd-order $G$, this holds for arbitrary $V'$ when
		$p\notin S(G)$.
	\end{corollary}
	\begin{proof}
		For each $g$, skew-primitive representatives of $(V_H)_g$ may
		be chosen in $P_{g,1}(H)$. Since $G$ is abelian, conjugation
		preserves this space. Averaging a linear section over $G$ makes
		the lifts equivariant, because $p\nmid|G|$.
		
		Give the chosen skew-primitive lifts filtration degree one and
		$\kk G$ filtration degree zero. The resulting algebra filtration
		$F_\bullet H'$ is a Hopf algebra filtration: the skew-primitive
		coproducts and antipodes of these lifts preserve it. Its associated
		graded algebra has a degree-zero Hopf projection, and
		\[
		\gr_FH'\simeq S\#\kk G.
		\]
		The right $G$-translates of these lifts have images in $H_1/H_0$
		spanning $V'\otimes\kk G$, so $F_1H'/F_0H'\simeq V'\otimes\kk G$
		and $S(1)=V'$.
		Since $\gr_FH'$ is generated in degrees zero and one,
		$T(V')\twoheadrightarrow S$ is a braided Hopf algebra epimorphism.
		Thus $S$ is pre-Nichols.
		
		The Nichols--Zoeller theorem gives $\dim H'\mid\dim H$, and hence
		\[
		\dim S=\frac{\dim H'}{|G|},\qquad p\nmid\dim S.
		\]
		Theorem~\ref{thm:uniform} gives $S=\B(V')$.
		For odd-order $G$, use Corollary~\ref{cor:odd-uniform}.
	\end{proof}
	
	This is the positive-characteristic dimension identity
	corresponding to \cite[Lemma~3.6]{AA} under the stated restrictions
	on the characteristic.
	
	For a Hopf algebra $C$ and a normalized convolution-invertible
	Hopf $2$-cocycle $\sigma:C\otimes C\to\kk$, write $C_\sigma$
	for the Hopf algebra with the same coalgebra and multiplication
	\[
	a\cdot_\sigma b
	=\sigma(a_1,b_1)a_2b_2\sigma^{-1}(a_3,b_3).
	\]

	\begin{corollary}\label{cor:cocycle-deformation}
		For every finite abelian group $G$ of odd order there is a finite
		set of rational primes $S_{\mathrm{def}}(G)\supseteq
		S_{\mathrm{odd}}(G)$ with the following
		property. Let $\kk$ be algebraically closed of characteristic
		$p\notin S_{\mathrm{def}}(G)$, and let $H$ be a finite-dimensional
		pointed Hopf algebra such that $G(H)\simeq G$ and $p\nmid\dim H$.
		Then, for some normalized convolution-invertible Hopf $2$-cocycle
		$\sigma$,
		\[
		H\simeq(\gr H)_\sigma
		\simeq\bigl(\B(V_H)\#\kk G\bigr)_\sigma.
		\]
	\end{corollary}
	\begin{proof}
		Put $d_G=2|G|^2$. By Corollary~\ref{cor:pointed},
		$\gr H\simeq\B(V_H)\#\kk G$ outside $S(G)$.
		The proof of Corollary~\ref{cor:odd-uniform} and
		Proposition~\ref{prop:transfer} give $\dim V_H\leq d_G$ and a
		characteristic-zero lift of the weights of $V_H$ in the finite list
		$\mathcal V(G,d_G)$. For each member of this list, choose a model
		$B_A$ as in Proposition~\ref{prop:general-model} and apply
		Lemma~\ref{lem:cocycle-specialization} to $C_A=B_A\#AG$.
		By \eqref{eq:structural-primes}, $S(G)\subseteq S_{\mathrm{odd}}(G)$.
		Define $S_{\mathrm{def}}(G)$ as the union of $S_{\mathrm{odd}}(G)$
		and the rational primes excluded in these finitely many
		choices and applications.
		
		Every $V_H$ under consideration is the reduction of one of
		these weight data: reduction identifies the roots of unity of
		order dividing the exponent of $G$, since $p\nmid|G|$.
		The corresponding fiber $B_\kk$ is $\B(V_H)$ by
		Proposition~\ref{prop:general-model}; its defining rank conditions
		remain valid after extending the residue field to $\kk$.
		Thus $H$ is a lifting of $C_\kk$, and the lemma applies.
		The rigidity exceptional set and the finitely many exclusions
		from the integral forms and specialization depend only on $G$,
		so $S_{\mathrm{def}}(G)$ does as well.
	\end{proof}
	
	\raggedbottom
	\appendix
	\section{Weights of the defining relations}\label{app:weights}
	
	Retain the characteristic-zero weight data of
	Section~\ref{sec:diagonal}. We verify
	\eqref{eq:relation-weight-exclusion} for homogeneous generators
	of $J(V)$. The relations in \cite[Theorem~3.1]{Angiono} span
	$\Rel(V)$. In each multidegree choose a subset whose classes
	form a basis, and denote the resulting set by $\mathcal R$.
	Induction on total degree shows that $\mathcal R$ generates
	$J(V)$: a relation with zero minimal class is a sum of products
	involving relations of smaller degree.
	
	For $\rho\in\mathcal R$ of multidegree
	$\boldsymbol a=(a_1,\ldots,a_\theta)$, write
	\[
	g_\rho=\prod_i g_i^{a_i},\qquad
	\chi_\rho=\prod_i\chi_i^{a_i},\qquad
	\gamma_\rho=(g_\rho,\chi_\rho).
	\]
	The coefficients involving the added weight in
	$V\oplus K_{\gamma_\rho}$ are
	\[
	q_{i\rho}=\prod_jq_{ij}^{a_j},\qquad
	q_{\rho i}=\prod_jq_{ji}^{a_j},\qquad
	q_{\rho\rho}=\prod_{i,j}q_{ij}^{a_i a_j}.
	\]
	If $q_{\rho\rho}=1$, the added line has Nichols algebra $K[t]$,
	since $\operatorname{char}K=0$. In particular its Nichols algebra,
	and hence that of the direct sum, is infinite-dimensional.
	For \eqref{eq:relation-weight-exclusion} it remains to consider
	$q_{\rho\rho}\ne1$.
	
	A nonzero minimal relation class satisfies
	\begin{equation}\label{eq:relation-antipode-weight}
		q_{\rho\rho}=\prod_iq_{ii}^{a_i}.
	\end{equation}
	To prove this, set $n=\sum_i a_i$ and use the braided Hopf quotient
	$A_n=T(V)/\langle J(V)(j):j<n\rangle$.
	The class of $\rho$ in $A_n$ is nonzero and primitive by
	\eqref{eq:first-relation-object}. Lemma~\ref{lem:antipode-twist}
	makes $\Theta^{-1}S^2$ an algebra endomorphism. Its value on
	$x_i$ is $q_{ii}^{-1}x_i$, so on multidegree $\boldsymbol a$
	it acts by $\prod_iq_{ii}^{-a_i}$. On the primitive class of
	$\rho$ it acts by $q_{\rho\rho}^{-1}$, which proves the formula.
	We also use $J(V)\cap T(U)=J(U)$ for a subspace $U$ spanned by
	part of the weight basis. Nonnegative multigrading identifies
	the corresponding minimal relation components in degrees
	supported on $U$. Thus redundancy can be checked on the
	generators occurring in a relation.
	
	For the root-vector powers (40), (42) and (43) of
	\cite[Theorem~3.1]{Angiono},
	\[
	\rho=x_\beta^{N_\beta},\qquad
	q_{\rho\rho}=q_{\beta\beta}^{N_\beta^2}=1.
	\]
	If $U$ is spanned by part of the weight basis of
	$V\oplus K_{\gamma_\rho}$, then $\B(U)$ embeds in
	$\B(V\oplus K_{\gamma_\rho})$. Finite-dimensionality of the latter
	would make the PBW root set of $\B(U)$ finite. These roots define
	an arithmetic root system by
	\cite[Theorem~1, arXiv version]{HeckRank2}.
	We apply the classification to these braided subspaces.
	For the remaining relations we use the following correspondence.
	The entries refer to Angiono's numbering.
	\begin{center}
		\small
		\begin{tabular}{ll}
			\toprule
			Relations in \cite[Theorem~3.1]{Angiono}
			& Results in \cite{Angiono}\\
			\midrule
			(41) & Proposition 4.1\\
			(44) & Lemma 4.2\\
			(45), (46) & Lemma 4.3\\
			(47) & Lemma 4.4\\
			(48), (49) & Lemma 4.5\\
			(56) & Lemma 4.6\\
			(50), (53) & Lemma 4.7\\
			(57), (58) & Lemma 4.8\\
			(59) & Lemma 4.9\\
			(51), (52) & Lemma 4.10\\
			(54) & Lemma 4.11(i)\\
			(55), (60), (61) & Lemma 4.11(ii)--(iv)\\
			(62) & Lemma 4.12(i)\\
			(63) & The calculation below\\
			(64), (65) & Lemma 4.12(iii), (iv)\\
			(66)--(68) & Lemma 4.12(v)--(vii)\\
			\bottomrule
		\end{tabular}
	\end{center}
	The cited proofs calculate the generalized Dynkin diagram
	obtained by adjoining a vector of the relation's weight and
	show that the resulting Nichols algebra is infinite-dimensional. These calculations depend
	only on the weights and therefore apply to the corresponding
	braided subspace of $V\oplus K_{\gamma_\rho}$.
	For (41), the condition $q_{ii}^{m_{ij}+1}\ne1$ is precisely
	that of \cite[Proposition~4.1]{Angiono}.
	For rank two and rank three we use
	\cite[Section~3, Tables~1 and~2]{HeckArithmetic}, up to
	permuting the vertices, and compute $m_{ij}=-a_{ij}$ by
	\eqref{eq:cartan}.
	Set $\widetilde q_{ij}=q_{ij}q_{ji}$.
	
	If $\widetilde q_{ij}=1$ in (44), then also
	$\widetilde q_{jk}=1$; the relation lies in the ideal generated
	by the degree-two commutators. Otherwise the hypotheses of
	\cite[Lemma~4.2]{Angiono} hold.
	
	In the fourth parameter case of (48), put $a=q_{ii}$ and
	$d=q_{kk}$. Here $q_{jj}=-1$,
	$\widetilde q_{ij}=a^{-2}$,
	$\widetilde q_{jk}=-a^{-3}$ and $\widetilde q_{ik}=1$.
	A nonzero minimal class of degree
	$2\alpha_i+3\alpha_j+\alpha_k$ would satisfy
	\[
	a^{-17}d=-a^2d,\qquad a^{19}=-1
	\]
	by \eqref{eq:relation-antipode-weight}.
	If $a=-1$, the degree-two relation $x_{ij}$ makes this class zero.
	Otherwise $\ord(a)=38$. In the basis obtained by reflection at $j$,
	the $i,k$ vertices define a rank-two arithmetic root subsystem
	\cite[Proposition~6]{HeckArithmetic}, with
	\[
	q'_{ii}=-a^{-1},\qquad
	\widetilde q'_{ik}=-a^{-5}=(q'_{ii})^5,\qquad
	\ord(q'_{ii})=19.
	\]
	No connected rank-two diagram in \cite[Table~1]{HeckArithmetic}
	has these parameters: for a vertex of order $19$, its edge label
	must be its inverse, inverse square or inverse cube.
	This rules out a nonzero minimal class in this parameter case.
	The other cases of (48), and (49), satisfy the hypotheses of
	\cite[Lemma~4.5]{Angiono}.
	
	For (50) and (53), equation~\eqref{eq:relation-antipode-weight}
	gives respectively
	$\widetilde q_{ij}^{\,3}=1$ and
	$\widetilde q_{ij}^{\,4}=1$.
	Together with their parameter conditions, these equalities give
	$q_{jj}=\widetilde q_{ij}^{-1}=\widetilde q_{jk}$,
	as required in \cite[Lemma~4.7]{Angiono}.
	
	For (54), put $b=q_{jj}$ and $d=q_{kk}$. Its diagram on
	$x_i,x_j,x_k$ is a path with vertex labels $-1,b,d$ and
	edge labels $-1,b^{-1}$, where $b\ne-1$.
	The complete rank-three list
	\cite[Table~2]{HeckArithmetic} gives exactly the following
	possibilities: row~3 gives $d=b$, $\ord(b)=4$;
	rows~6, 7 and~17 give $d=-1$, with
	$\ord(b)=4,6,3$, respectively. This comparison includes both
	orientations of every path in the table.
	When $d=-1$, \cite[Lemma~4.11(i)]{Angiono} applies.
	When $d=b$ has order four, the diagram is of Cartan type $C_3$
	on the ordered generators $x_k,x_j,x_i$.
	The complete presentation \cite[Section~4.3.5]{AADiagonal}
	has no generator of multidegree
	$2\alpha_i+2\alpha_j+\alpha_k$; its minimal relation component
	in that degree is therefore zero.
	
	For (62)--(65), write $a=q_{ii}$, $b=q_{jj}$ and
	$t=\widetilde q_{ij}$.
	If $t=1$, relation (62) is a multiple of $x_{ij}^2$ and has
	zero minimal class, since $x_{ij}\in J(V)(2)$.
	If $t\ne1$, division by $(1-t)bq_{ji}$ gives the relation in
	\cite[Lemma~4.12(i)]{Angiono}.
	
	Relation (63) has multidegree $4\alpha_i+2\alpha_j$.
	For its branch $m_{ij}\in\{4,5\}$, the edge labels to the added
	vertex are
	\[
	u=q_{i\rho}q_{\rho i}=a^8t^2,\qquad
	v=q_{j\rho}q_{\rho j}=t^4b^4,\qquad
	q_{\rho\rho}=u^2v.
	\]
	If $q_{\rho\rho}\ne1$, at least one of $u,v$ differs from $1$.
	The enlarged rank-three diagram is connected and contains the
	rank-two diagram on $x_i,x_j$, with $m_{ij}\ge4$.
	Every rank-two subdiagram in \cite[Table~2]{HeckArithmetic}
	has $m_{ij}\le3$, excluding this case.
	In the other branch, $m_{ij}=3$, $\ord(a)=4$ and $b=-1$,
	equation~\eqref{eq:relation-antipode-weight} gives
	$q_{\rho\rho}=a^4b^2=1$ for a nonzero minimal class.
	
	For (64), the hypotheses with $m_{ji}>2$ leave only
	\[
	(a,t,b)=(\zeta^8,-\zeta^{-1},\zeta^6),\quad \ord(\zeta)=24,
	\qquad
	(a,t,b)=(-\eta^5,\eta^{-3},\eta),\quad \ord(\eta)=30,
	\]
	in rows~12 and~15 of \cite[Table~1]{HeckArithmetic}, with table
	parameters $\zeta$ and $-\eta$, respectively.
	They are obtained by interchanging the vertices in the diagrams
	\cite[(10.56 a), (10.64 d)]{AADiagonal}.
	After interchanging the generators, the complete presentations
	\cite[(10.57), (10.68)]{AADiagonal} have no generator of
	multidegree $4\alpha_i+3\alpha_j$.
	Thus this component of $\Rel(V)$ is zero.
	The remaining parameters have $b=-1$ or $m_{ji}=2$ and are
	covered by \cite[Lemma~4.12(iii)]{Angiono}.
	
	For (65), the degree is $5\alpha_i+3\alpha_j$, so a nonzero
	minimal class requires
	\[
	a^{20}b^6t^{15}=1.
	\]
	The proof of \cite[Lemma~4.12(iv)]{Angiono} reduces the possible
	nonzero minimal classes to the parameters
	\[
	(a,t,b)=(\zeta^3,\zeta^8,-1),\quad \ord(\zeta)=9,
	\qquad
	(a,t,b)=(\eta^3,-\eta^4,-\eta^{-4}),\quad \ord(\eta)=15.
	\]
	In the latter case,
	$q_{\rho\rho}=\eta^9\ne-\eta^3=a^5b^3$, so the minimal class
	is zero.
	In the order-nine case the class is nonzero.
	The complete presentation \cite[(8.3)]{AADiagonal}
	has only $x_i^3$ and $x_j^2$ in total degrees less than eight.
	The relation $[x_{iij},[x_{iij},x_{ij}]_c]_c$ has coefficient
	$1$ on $(x_i^2x_j)^2x_ix_j$, which contains neither of these
	words, so it does not belong to the ideal they generate.
	Here
	\[
	q_{\rho\rho}=-\zeta^6,\qquad
	q_{j\rho}q_{\rho j}=\zeta^4\ne1.
	\]
	The diagram on $x_i,x_j$ occurs in
	\cite[Table~1, row~9]{HeckArithmetic}, but in none of the connected
	rank-three diagrams in \cite[Table~2]{HeckArithmetic}.
	Since $q_{j\rho}q_{\rho j}\ne1$, adjoining the relation weight
	gives a connected diagram, so the enlarged Nichols algebra is
	infinite-dimensional.
	
	These checks, together with the diagram exclusions in the proofs
	cited in the correspondence table, establish \eqref{eq:relation-weight-exclusion}
	for every member of $\mathcal R$.
	
	\section{Specialization of cocycle deformations}
	\label{app:cocycle}
	
	\begin{lemma}\label{lem:cocycle-specialization}
		Let $G$ be a finite abelian group, $K$ a number field, and
		$A$ a localization of its ring of integers such that $|G|$
		is invertible in $A$. Suppose
		\[
		B_A=\bigoplus_{n=0}^h B_A(n)
		\]
		is a finite free connected graded braided Hopf algebra in
		${}_{AG}^{AG}\mathcal{YD}$ whose geometric fibers are Nichols
		algebras, and put $C_A=B_A\#AG$.
		
		There is a finite set of rational primes $S$ such that, for
		every algebraically closed field $\kk$ equipped with a map
		$A\to\kk$ and with $\operatorname{char}\kk\notin S$, every
		finite-dimensional filtered Hopf algebra $L$ satisfying
		\[
		\gr_F L\simeq C_\kk
		\]
		as graded Hopf algebras is isomorphic to $(C_\kk)_\sigma$
		for some Hopf $2$-cocycle $\sigma$.
	\end{lemma}
	
	\begin{proof}
		After inverting a nonzero element of $A$, each of the finitely many
		projective graded components $C_A(n)$ is free. Include the rational
		primes below the removed prime ideals in $S$; every map $A\to\kk$
		in a remaining characteristic factors through this localization.
		Choose a homogeneous basis of $C_A$ over the localized ring, still
		denoted by $A$, and put
		\[
		F_nC_A=\bigoplus_{j\le n}C_A(j).
		\]
		On the underlying free module, consider Hopf algebra
		structures preserving this filtration, with the prescribed
		Hopf structure on $F_0=AG$ and associated graded Hopf
		algebra $C_A$. In the chosen basis, filtration preservation,
		the Hopf axioms, and the prescribed associated graded
		structure are finitely many polynomial equations in the
		structure constants. They define an affine scheme $X_A$
		of finite type. Every $L$ in the statement determines a
		$\kk$-point of $X_A$ after choosing compatible lifts of a
		graded basis.
		
		For a geometric point of $X_A$, let $L$ denote the corresponding
		filtered Hopf algebra and put
		\[
		J=F_0^\perp\subseteq L^*.
		\]
		Since $F$ is a coalgebra filtration,
		\[
		J^{h+1}=0,\qquad
		L^*/J\simeq(\kk G)^*\simeq\kk^G.
		\]
		Thus $J=\operatorname{rad}(L^*)$, and hence $L_0=F_0$.
		The coalgebra filtration gives $F_n\subseteq L_n$.
		
		Conversely, suppose $x\in L_n$ has filtration degree
		$d>n$, and let $\bar x$ be its nonzero image in
		$(\gr_F L)(d)$. By the iterated coproduct criterion
		for the coradical filtration,
		\[
		(\pi^{\otimes(n+1)})\Delta^{(n)}(x)=0,
		\qquad \pi:L\to L/L_0.
		\]
		Passing to the associated graded shows that $\bar x$
		belongs to the $n$th coradical term of $\gr_F L$.
		But $\gr_F L\simeq C_\kk$ is coradically graded,
		so this term is
		$\bigoplus_{j\le n}C_\kk(j)$, contradicting $d>n$.
		Therefore
		\[
		F_n=L_n\qquad(n\ge0).
		\]
		
		Next adjoin to the coordinates of $X_A$ the coefficients
		of a linear map
		\[
		\sigma:C_A\otimes_A C_A\to A,
		\]
		its convolution inverse, and linear maps in both directions
		between $(C_A)_\sigma$ and the corresponding Hopf algebra.
		The normalization and cocycle identities, the two convolution
		inverse identities, and the equations requiring the linear maps
		to be mutually inverse and to preserve multiplication, coproduct,
		unit and counit are finitely many polynomial equations.
		In particular, the multiplication
		of $(C_A)_\sigma$ is polynomial in these coordinates.
		We obtain a finite-type morphism
		\[
		f:Y_A\longrightarrow X_A.
		\]
		
		Over $K=\operatorname{Frac}(A)$, every geometric point
		of $X_K$ represents a finite-dimensional pointed Hopf
		algebra with abelian group of group-like elements and
		associated graded Hopf algebra isomorphic to a geometric fiber of $C_A$. Since the geometric
		fibers of $B_A$ are Nichols algebras,
		\cite[Corollary~1.2]{AGI} shows that every such Hopf
		algebra is a cocycle deformation of that fiber of $C_A$.
		Thus $f_K$ is surjective on the underlying topological
		spaces.
		
		By Chevalley's theorem (see e.g. \cite[Tag~054K]{Stacks}),
		\[
		E=X_A\setminus f(Y_A)
		\]
		is constructible, as is its image under
		$X_A\to\operatorname{Spec}A$. This image does not contain
		the generic point of $\operatorname{Spec}A$. Since $A$
		is a localization of the ring of integers of a number
		field, the image is finite. After excluding the rational
		primes below these closed points, the morphism $f$
		is surjective.
		
		For every remaining algebraically closed field $\kk$
		as in the statement, the fiber of $f$ over any
		$\kk$-point of $X_A$ is a nonempty scheme of finite
		type over $\kk$, and hence has a $\kk$-point. Its
		coordinates give a Hopf $2$-cocycle $\sigma$ and an
		isomorphism $L\simeq(C_\kk)_\sigma$.
	\end{proof}

\end{document}